\documentclass[a4paper,12pt]{amsart}
\newcommand{\topcounter}{section}

\usepackage[francais]{babel}
\usepackage[applemac]{inputenc}

\usepackage[a4paper]{geometry}
\usepackage{graphicx}
\usepackage{amsmath}
\usepackage{amsthm}
\usepackage{amssymb}
\usepackage{amsfonts}
\usepackage{mathrsfs}
\usepackage{enumerate}
\usepackage{fancyhdr}

\theoremstyle{plain}
\newtheorem{theoreme}{Thorme}[\topcounter]

\newtheorem{dictionnaire*}{Dictionnaire}

\newtheorem{lemme}[theoreme]{Lemme}
\newtheorem*{lemme*}{Lemme}

\newtheorem*{corollaire*}{Corollaire}
\newtheorem*{important}{Important}

\theoremstyle{definition}
\newtheorem{definition}{Dfinition}[\topcounter]

\newtheorem{principe}{Principe}

\newtheorem*{demonstration}{Dmonstration}
\newtheorem{demonstration*}{Dmonstration}

\newtheorem{exemple}{Exemple}
\newtheorem*{exemples}{Exemples}

\newtheorem*{convention*}{Convention}

\newtheorem*{generalisation*}{Gnralisation}
\newtheorem{application}{Application}
\newtheorem*{application*}{Application}

\theoremstyle{remark}

\newtheorem*{remarque}{Remarque}
\newtheorem*{remarquesubsidiaire}{Remarque subsidiaire}

\newtheorem*{variante*}{Variante}

\providecommand{\abs}[1]{\lvert \, #1 \, \rvert}

\def\inv{^{-1}}

\DeclareMathOperator{\id}{id}

\DeclareMathOperator{\image}{Im}

\newcommand{\idest}{\textit{i.e.} }

\newcommand{\confer}{\textit{cf.} }

\newcommand{\ZZ}{\mathbb{Z}}

\newcommand{\CC}{\mathbb{C}}

\newcommand{\Gr}{\Gamma}

\newcommand{\shortpar}{\vspace{-4mm}}

\renewcommand{\emptyset}{\varnothing}

\renewcommand{\epsilon}{\varepsilon}

\input xy
\xyoption{all}

\newcommand{\ZSeries}{\mathcal{Z}}
\newcommand{\bary}{^\mathrm{sb}}
\newcommand{\baryplus}{^{\mathrm{sb}+}}

\newcommand{\cyclique}{{\ZZ}}
\newcommand{\cycliquedeux}{{\ZZ/2\ZZ}}
\newcommand{\graphique}{{\cyclique\ast\cycliquedeux}}

\renewcommand{\phi}{\varphi}
\renewcommand{\Im}{\image}

\DeclareMathOperator{\Fix}{Fix}
\DeclareMathOperator{\Sym}{Sym}

\subjclass[2000]{Primary 05C25, 05C30, 05C85, 20F05 ; Secondary 20F10, 05C38, 20F36}

\begin{document}

\title{Sur la Classification et le Denombrement des Sous-groupes du Groupe Modulaire et de leurs Classes de Conjugaison}
\author{Samuel Alexandre \textsc{Vidal}}

\maketitle

\renewcommand{\abstractname}{Abstract}

\begin{abstract}

In this article we give a  classification of the sub-groups in $\mathrm{PSL}_2(\ZZ)$ and of the conjugacy classes of these sub-groups by the mean of an combinatorial invariant~: some trivalent diagrams (dotted or not). We give explicit formulae enabling to count the number of isomorphism classes of these structures and of some of their variations, as function of the number of their arcs. Until now, the counting of non-dotted diagrams was an open problem, for it gives also the number of unrooted combinatorial maps, triangular or general respectively. The article ends with the description of a high performance algorithm to enumerate those structures witch is built upon an unexpected factoring of the cycle index series of the considered combinatorial species.

\end{abstract}

\renewcommand{\abstractname}{Rsum}

\begin{abstract}
Dans cet article nous donnons une classification des sous-groupes de $\mathrm{PSL}_2(\ZZ)$ et les classes de conjugaison de ces sous-groupes au moyen d'un invariant de nature combinatoire~: des diagrammes trivalents (points ou non). Nous donnons des formules explicites permettant de compter les classes d'isomorphisme de ces structures et de plusieurs de leurs variantes, en fonction du nombre de leurs arcs. Jusqu' prsent, le dnombrement de diagrammes non-points tait un problme ouvert du fait que cela donne le nombre de cartes combinatoire non-enracines, triangulaires ou quelconques respectivement. L'article s'achve par la description d'un algorithme extrmement efficace pour numrer ces structures qui s'appuie sur une factorisation inattendue de la srie indicatrice des cycles des espces combinatoires considres.
\end{abstract}

\setcounter{section}{-1}


\section{Introduction}

Les sous-groupes de $\mathrm{PSL}_2(\ZZ)$ interviennent par exemple dans l'tude des quations fonctionnelles de fonctions multivalues sur des courbes algbriques projectives \idest surfaces de Riemann compactes, notemment dans l'tude des fonctions Polylogarithmes gnralises sur des revtements algbriques de la sphre de Riemann. On sait par ailleurs qu'elles correspondent de faon exhaustive \cite{belyi80}, aux courbes algbriques projectives (munies d'un nombre fini de \emph{points marqus}) dfinies sur la clture algbrique $\bar{\mathbb{Q}}$ du corps des nombres rationnels $\mathbb{Q}$. L'tude de ces structures en vue d'une meilleure comprhension du groupe de Galois abolu $\mathrm{Gal}_\mathbb{Q}(\bar{\mathbb{Q}})$ constitue un des points les plus saillants du vaste programme galoisien d'A. Grothendieck \cite{grothendieckprogramme} comme en tmoigne le nombre de publications de haut niveau consacres  ces questions sur les vingt dernires annes. On peut mme faire remonter l'origine de ce sujet aux traveaux de F. Klein sur la gomtrie de l'icosadre \cite{klein56}.

Le point principal se trouve sans doute tre que $\mathrm{Gal}_\mathbb{Q}(\bar{\mathbb{Q}})$ agit \emph{fidlement} sur la catgorie de ces revtements (en modifiant les donnes de monodromie). Cette action est toutefois encore trs loin d'tre bien comprise et cela en dpit des nombreux efforts qui ont t consacrs  son tude. Elle jouit d'excellentes proprits de finitude, ce qui attise l'espoir de parvenir un jour  la dcrire compltement. Les objets permuts sont en effet des structures combinatoires basiques et les orbites de cette action sont toutes finies. Un nombre important de publications ont t consacres et sont encore consacres  l'tude d'invariants prservs par cette action de groupe dans l'espoir d'en caractriser un jour les orbites avec prcision.

L'tude des formes modulaires s'est pour l'instant beaucoup concentre sur la famille des sous-groupes de congruences prsents dans $\mathrm{PSL}_2(\ZZ)$ lesquels reoivent une interprtation modulaire astucieuse de  P. Deligne \cite{deligne73}. Ils sont associs  la catgorie des courbes elliptiques classifies modulo leurs isognies centralises par leurs points de $n$-division.

\begin{figure}[h]
\begin{center}
\includegraphics{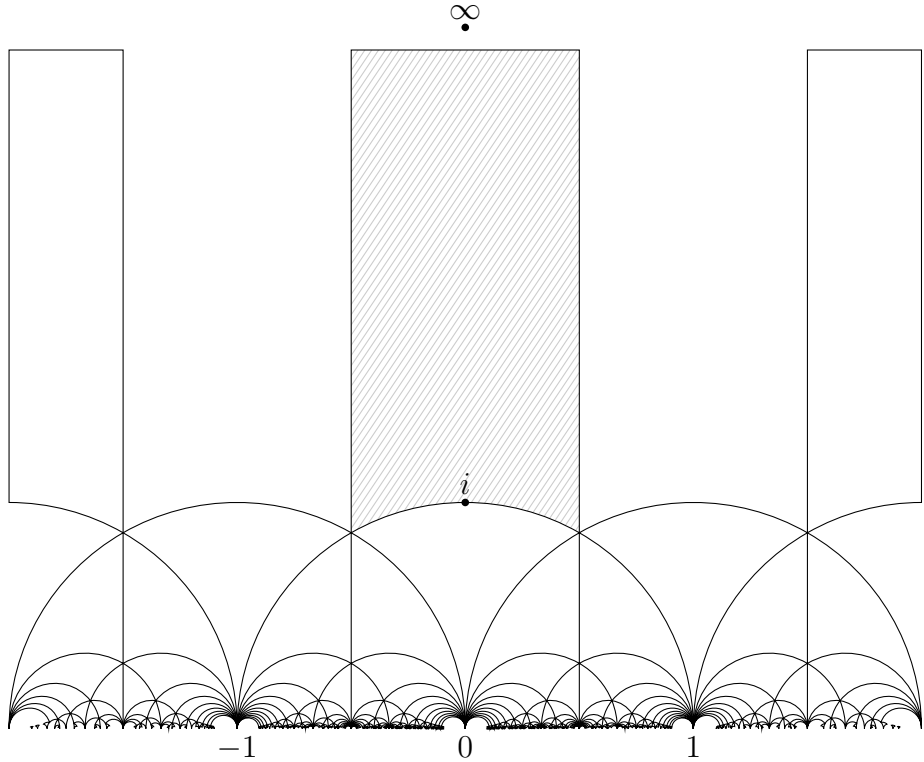}
\caption{Pavage du demi-plan de Poincar par les translats du domaine fondamental $\mathscr{F}$ sous l'action du groupe modulaire $\mathrm{PSL}_2(\ZZ)$.}
\label{groupe:modulair:pavage}
\end{center}
\end{figure}

Le groupe modulaire se conoit galement comme l'exemple le plus significatif de groupe Fuschien ainsi qu'en tmoignent les illustres travaux d'H. Poincar. Son action sur le demi-plan de Poincar $\mathscr{H}$ est donne par les transformation homographiques
\begin{align*}
	&& \pm \begin{pmatrix} a & b \\ c & d \end{pmatrix}\cdot \tau =
		\frac{a\,\tau+b}{c\,\tau+d}
	&& \Im \tau > 0
\end{align*}
avec $a$, $b$, $c$ et $d$ entiers rationnels vrifiant $ad- bc = 1$. Cette action de groupe admet un domaine fondamental $\mathscr{F}$ qui produit un pavage du demi-plan de Poincar par des triangles hyperboliques.
La figure \ref{groupe:modulair:pavage} reprend cette situation qui est des plus classiques.

Les lments $\tau$ du demi-plan de Poincar $\mathscr{H}$ sont senss reprsenter l'invariant $\omega_2/\omega_1$, o $(\omega_1,\omega_2)$ dcrit une base oriente d'un rseau de $\CC$. C'est un invariant complet pour l'action du groupe  $\CC^*$ des homothties complexes. Une telle base est dite \emph{rduite} lorsqu'elle satisfait aux conditions de Gauss, qui exprimes en terme de l'invariant $\tau$ se rduisent aux ingalits suivantes~:
\begin{align*}
	\abs{1} &\le \abs{\tau}				&
	\abs{\tau - 1} &\le \abs{\tau}		&
	\abs{\tau + 1} &\le \abs{\tau}
\end{align*}
et ce sont prcisment celles-ci qui dfinissent le domaine fondamental $\mathscr{F}$ de l'action du groupe modulaire $\mathrm{PSL}_2(\ZZ)$ sur le demi-plan de Poincar.

Dans ce travail, nous parvenons  classifier les sous-groupes du groupe modulaire $\mathrm{PSL}_2(\ZZ)$, ainsi que leurs classes de conjugaison, au moyen d'invariants de nature combinatoire : des diagrammes. Notre argument s'applique en principe, moyennant quelques adaptations triviales, au cas du produit libre de deux groupes cycliques quelconques. Par exemple, le groupe modulaire $\mathrm{PSL}_2(\ZZ)$ est produit libre des deux groupes cycliques $\ZZ/2\ZZ$ et $\ZZ/3\ZZ$ et ses sous-groupes se classifient de fait au moyen de diagrammes trivalents.
On donne des critres simples sur les diagrammes permettant de dcider rcursivement (\idest au moyen d'un algorithme) si un sous-groupe est distingu, si deux sous-groupes sont conjugus, ou plus simplement si deux sous-groupes sont en relation d'inclusion.
On parvient en particulier aux correspondances biunivoques suivantes :
\begin{center}
\begin{tabular}{ccccc}
	\begin{tabular}{c}
	Diagrammes trivalents \\  isomorphismes prs
	\end{tabular}
	& $\longleftrightarrow$ &
	\begin{tabular}{c}
	Sous-groupes de $\mathrm{PSL}_2(\ZZ)$ \\  conjugaison prs
	\end{tabular}
	\\ \\
	\begin{tabular}{c}
	Diagrammes trivalents points \\  isomorphismes prs
	\end{tabular}
	& $\longleftrightarrow$ &
	\begin{tabular}{c}
	Sous-groupes de $\mathrm{PSL}_2(\ZZ)$
	\end{tabular}
\end{tabular}
\end{center}

Fort de ces quivalences, nous obtenons des formules gnrales, sous forme de sries gnratrices pour compter le nombre de sous-groupes d'indice fini donn dans $\mathrm{PSL}_2(\ZZ)$ et le nombre de leurs classes de conjugaison. Elles concident en effet avec le nombre de digrammes trivalents et leurs variantes points pour un nombre d'arcs donn. On aboutit dans un premier temps  la formule gnrale donnant le nombre de diagrammes trivalents \emph{points}, laquelle conduit aux premires valeurs suivantes,
\begin{align*}
	\tilde{D}_3^\bullet(t) = t+{t}^{2}+4\,{t}^{3}+8\,{t}^{4}+5\,{t}^{5}+22\,{t}^{6}+42\,{t}^{7}
	+40\,{t}^{8}+120\,{t}^{9}
	+\dots
\end{align*}
Lesquelles avaient dj calcules par  Stothers en 1977 \cite{stothers77}.

La formule gnrale donnant le nombre de diagrammes trivalents \emph{non-points} de taille donne laquelle permet de compter le nombre de \emph{classes de conjugaison} de sous-groupes ayant un indice fini donn dans $\mathrm{PSL}_2(\ZZ)$,
\begin{align*}
	\tilde{D}_3(t)\, &= \,t+{t}^{2}+2\,{t}^{3}+2\,{t}^{4}+{t}^{5}+8\,{t}^{6}+6\,{t}^{7}+7\,{t}^{8}+14\,{t}^{9}
	+\dots 
\end{align*}
tait reste jusqu' maintenant un problme ouvert et elle est nettement plus difficile  obtenir. La solution que nous prsentons, fait intervenir de faon cruciale la thorie des \emph{espces combinatoires}, due  A. Joyal et aux travaux de l'cole qubcoise de combinatoire.
Les tables \ref{tab:diag:3:taille:inf:5}  \ref{tab:diag:3:taille:9} pages  \pageref{tab:diag:3:taille:inf:5}  \pageref{tab:diag:3:taille:9} donnent une liste exhaustive des diagrammes trivalents \emph{non points} de taille infrieure  \emph{neuf}. Elles vrifient prcisment ce dcompte.

On dispose galement d'une interprtation des diagrammes trivalents en terme de revtements et la reprsentation de monodromie correspondante permet une relecture de ce dnombrement en terme de permutations~: le $n$-ime terme correspond au nombre de paires de permutations $(\tau_1,\tau_2)$  conjugaison simultane prs qui agissent transitivement sur un ensemble de taille $n$ et qui vrifient les conditions d'\emph{involutivit} $\tau_1^2=\id$ et de \emph{triangularit} $\tau_2^3=\id$.

Enfin, la formule gnrale, crite navement, pose de srieux problmes de complexit~: les calculs ncessaires pour en dterminer les coefficients interdisent par leur volume de considrer des termes au del des tout premiers. Nous nous sommes donc attach dans la dernire partie de cet expos  donner une mthode de calcul efficace. Celle-ci s'appuie sur une factorisation inattendue de la srie gnratrice correspondante. Il en rsulte un effondrement de la complexit qui autorise un calcule en trs grand poids. On illustre ce point en donnant les termes de poids cinq-cent des sries considres.

\section{Prliminaires Combinatoires}

\subsection{Graphes Gnralits}

Dans cette section nous introduisons une variante de la notion de graphe, laquelle servira de point de dpart  la construction de notre invariant combinatoire. La dfinition prsente ici diffre lgrement de la notion habituelle de graphe mais le texte signalera clairement ces quelques diffrences. On s'attachera en outre, dans la section suivante  dcrire une opration permettant de construire un graphe au sens habituel  partir d'un graphe au sens de la dfinition ci-dessous, et la donne d'une structure supplmentaire sur le graphe ainsi construit permettra de retrouver l'quivalence avec notre dfinition.


\begin{definition}
\label{def:graphes}
Par un \emph{graphe} $\Gr$, on entend la donne de deux ensembles $\Gr_0$ et $\Gr_1$, et de trois applications $s,t : \Gr_1 \to \Gr_0$ et $.\inv : \Gr_1 \to \Gr_1$. Le tout vrifiant pour tout $a\in \Gr_1$ :
\begin{align*}
	(a\inv)\inv &= a &
	s(a\inv) &= t(a) &
	t(a\inv) &= s(a)
\end{align*}
\end{definition}

Les \emph{sommets} du graphe sont les lments de $\Gr_0$, ses \emph{arcs} (ou \emph{demi-artes}) sont les lments de $\Gr_1$ et les deux applications $s$ et $t$ font correspondre  toute arc $a$, son \emph{origine} $s(a)$ et sa \emph{destination} $t(a)$.
Comme l'application $.\inv$ est involutive, on a une action du groupe  deux lments sur $\Gr_1$. Les \emph{artes} (ou \emph{bi-arcs}) du graphe sont les orbites dans $\Gr_1$ de cette action de groupe. On note $\Gr_1^*$ cet ensemble et $\pi$ la projection $\Gr_1 \to \Gr_1^*$.

\begin{remarque}
La dfinition que l'on donne n'exclut pas qu'un mme arc ait mme origine et mme destination, on dit que c'est une \emph{boucle}, ni qu'une arte ne comporte qu'un arc, auquel cas on dit que l'arte est \emph{plie} \cite[p. 15]{grothendieckprogramme}. En thorie quantique des champs, les diagrammes de Feynman sont des structures qui admettent de telles artes plies, pour lesquelles une extrmit reste ``ouverte", ce sont les \emph{artes externes} du diagramme ou \emph{pattes externes}.
\end{remarque}

\begin{exemple}
\label{fig:venus:1}
On observe que le graphe suivant admet \emph{deux} sommets et \emph{cinq} artes dont trois sont plies. Il y a donc \emph{sept} arcs.
Le tableau prcise les valeurs des trois applications correspondantes. On a pos $\Gr_0 = \{\,s_1,s_2\,\}$ et $\Gr_1 = \{\,a_1,a_2,a_3,a_4,a_5,a_4\inv,a_5\inv\,\}$.
\begin{center}
\begin{tabular}{cc}
	\begin{tabular}{c}
		\includegraphics{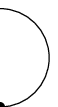}
	\end{tabular}
	&
	\begin{tabular}{l}
	\begin{tabular}{c|lllllll}
		& $a_1$ & $a_2$ & $a_3$ & $a_4$ & $a_5$ & $a_4\inv$ & $a_5\inv$	\\
	\hline
	$s$ & $s_2$ & $s_2$ & $s_2$ & $s_1$ & $s_1$ & $s_2$ & $s_1$	\\
	$t$ & $s_2$ & $s_2$ & $s_2$ & $s_2$ & $s_1$ & $s_1$ & $s_1$	\\
	$.\inv$ & $a_1$ & $a_2$ & $a_3$ & $a_4\inv$ & $a_5\inv$ & $a_4$ & $a_5$
	\end{tabular}
	\end{tabular}
\end{tabular}
\\
\end{center}
\end{exemple}

\begin{definition}
\label{def:graphes:morph}
Un \emph{morphisme} $\phi$ entre deux graphes $\Gr$ et $\Gr'$ est la donne de deux applications $\phi_0 : \Gr_0 \to \Gr'_0$ et $\phi_1 : \Gr_1\to \Gr'_1$ compatibles aux applications de structure en ce sens que les diagrammes suivants commutent :
\begin{align*}
\xymatrix@C=1.5cm@R=1.5cm
{
	{\Gr_1}
		\ar[r]^{\phi_1}
		\ar[d]_{s}			&
	{\Gr_1'}
		\ar[d]^{s}			\\
	{\Gr_0}
		\ar[r]_{\phi_0}		&
	{\Gr_0'}
}
&&
\xymatrix@C=1.5cm@R=1.5cm
{
	{\Gr_1}
		\ar[r]^{\phi_1}
		\ar[d]_{t}			&
	{\Gr_1'}
		\ar[d]^{t}			\\
	{\Gr_0}
		\ar[r]_{\phi_0}		&
	{\Gr_0'}
}
&&
\xymatrix@C=1.5cm@R=1.5cm
{
	{\Gr_1}
		\ar[r]^{\phi_1}
		\ar[d]_{.\inv}	&
	{\Gr_1'}
		\ar[d]^{.\inv}	\\
	{\Gr_1}
		\ar[r]_{\phi_1}		&
	{\Gr_1'}
}
\end{align*}
\end{definition}

\subsection{Subdivision Barycentrique}
\label{sec:subdivision:barycentrique}

Nous dcrivons maintenant un procd simple pour se dbarrasser des boucles et des artes plies. La \emph{subdivision barycentrique} est l'opration qui associe  tout graphe $\Gr$ le graphe $\Gr^\mathrm{sb}$ 
obtenu  partir de $\Gr$ en ajoutant un sommet supplmentaire au milieu de chaque arte.

De faon plus prcise, on associe  $\Gr$ le graphe $\Gr\bary$ avec $\Gr\bary_0 = \Gr_0 \sqcup \Gr_1^*$ et $\Gr\bary_1 = \Gr_1 \sqcup \Gr_1$. Les trois applications de structure sont dfinies par les trois diagrammes commutatifs suivants o l'on note $\rho_1$ et $\rho_2$ les injections naturelles associes  la runion disjointe $\Gr_1 \sqcup \Gr_1$ :
\begin{align*}
\xymatrix@C=.5cm@R=.7cm
{
	& {\Gr_1 \sqcup \Gr_1}
		\ar@{..>}[dd]^{s}		\\
	{\Gr_1}
		\ar[ur]^{\rho_1}
		\ar[dr]_{s}				&&
	{\Gr_1}
		\ar[ul]_{\rho_2}
		\ar[dl]^{\pi}			\\
	& {\Gr_0 \sqcup \Gr_1^*}
}
&&
\xymatrix@C=.5cm@R=.7cm
{
	& {\Gr_1 \sqcup \Gr_1}
		\ar@{..>}[dd]^{t}		\\
	{\Gr_1}
		\ar[ur]^{\rho_1}
		\ar[dr]_{\pi}			&&
	{\Gr_1}
		\ar[ul]_{\rho_2}
		\ar[dl]^{\text{$s$ (sic)}}			\\
	& {\Gr_0 \sqcup \Gr_1^*}
}
&&
\xymatrix@C=.5cm@R=.7cm
{
	& {\Gr_1 \sqcup \Gr_1}
		\ar@{..>}[dd]^{.\inv}		\\
	{\Gr_1}
		\ar[ur]^{\rho_1}
		\ar[dr]_{\rho_2}			&&
	{\Gr_1}
		\ar[ul]_{\rho_2}
		\ar[dl]^{\rho_1}			\\
	& {\Gr_1 \sqcup \Gr_1}
}
\end{align*}

L'opration est videmment \emph{fonctorielle} mais n'admet pas d'opration inverse, ne serait-ce que parce qu'une fois l'opration effectue, il n'y a plus moyen de distinguer les sommets qui proviennent de $\Gr_0$ de ceux qui proviennent de $\Gr_1^*$. Autrement dit, le foncteur n'est pas fidle. Une faon de remdier au problme consiste  enrichir la catgorie d'arrive par une information de \emph{couleur}.

\begin{definition}
Un \emph{graphe bicolori} est un graphe $\Gr$ muni d'une application $\alpha$ de l'ensemble de ses sommets  valeurs dans $\{\,0,1\,\}$.
\end{definition} \shortpar
Par convention, on qualifiera de \emph{blancs} (resp. de \emph{noirs}) les sommets $x$ tels que $\alpha(x)=0$ (resp. $\alpha(x)=1$). Les morphismes de graphes bicoloris sont les morphismes de graphes qui prservent cette information supplmentaire.

L'opration de \emph{subdivision barycentrique enrichie} est alors, l'opration qui  un graphe $\Gr$ associe le graphe $\Gr\bary$ muni de la coloration qui aux sommets provenant de $\Gr_0$ associe la valeur \emph{un} (noir) et  ceux qui proviennent de $\Gr_1^*$ associe la valeur \emph{zro} (blanc). On notera $\Gr\baryplus$ le graphe bicolori ainsi obtenu.
\begin{exemple}
On reprend dans le dessin de gauche le graphe de l'exemple \ref{fig:venus:1} plus haut. On reprsente au centre le rsultat de sa subdivision barycentrique et  droite le rsultat de sa subdivision barycentrique enrichie.
\begin{center}
\includegraphics{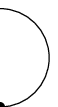}
\hspace{1.6cm}
\includegraphics{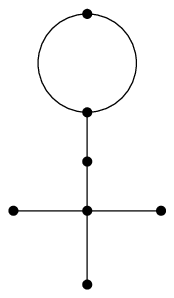}
\hspace{1.5cm}
\includegraphics{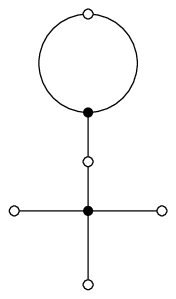}
\end{center}
\end{exemple}

Le foncteur $.\baryplus$ que l'on obtient par ce qui prcde est \emph{pleinement fidle}. Il n'est en revanche pas \emph{surjectif} pas mme au sens faible d'\emph{essentiellement surjectif}. Les sommets blancs sont en effet toujours de degr \emph{un} ou \emph{deux}, suivant qu'ils proviennent d'une arte plie ou non. De plus, deux sommets de mme couleur ne sont jamais relis par une mme arte.

\begin{definition}
Un graphe bicolori est dit \emph{propre} (\emph{clean}) si ses sommets blancs sont de degr \emph{un} ou \emph{deux} et si deux sommets de mme couleur ne sont jamais relis par une mme arte.
\end{definition}\shortpar

On vrifie sans peine que le foncteur de subdivision barycentrique enrichie
$.\baryplus$
est essentiellement surjectif sur la \emph{sous-catgorie pleine} des graphes bicoloris propres.
Il rsulte donc une quivalence de catgories entre cette dernire et la catgorie des graphes au sens des dfinitions \ref{def:graphes} et \ref{def:graphes:morph}.

\begin{remarquesubsidiaire}
Il y a une bijection naturelle vidente entre l'ensemble des arcs d'un graphe et les artes de sa subdivision barycentrique~; elle est d'ailleurs \emph{unique} ( conjugaison par $.\inv$ prs.)
\end{remarquesubsidiaire}

\subsection{Orientation Cyclique aux Sommets}

Dans ce qui suit, nous allons procder  quelques enrichissements suppl\-men\-taires de la catgorie des graphes.

\begin{definition}
L'\emph{toile} relative  un sommet $x$ d'un graphe $\Gr$ est par dfinition l'ensemble $s\inv(x) = \{\,a\in \Gr_1\,|\,s(a)=x\,\}$ des arcs ayant pour origine ce sommet.
\end{definition} \shortpar

Etant donn un graphe $\Gr$, une \emph{orientation cyclique aux sommets} est la donne d'une action $(a,k) \mapsto a+k$ du groupe cyclique infini $\cyclique$ sur l'ensemble des arcs de $\Gr$ qui premirement, soit compatible aux \emph{toiles} au sens o $s(a+1) = s(a)$ quel que soit l'arc $a$ et qui d'autre part, soit \emph{toile-transitive} au sens o quels que soient les arcs $a_1$ et $a_2$, la condition $s(a_1) = s(a_2)$ entrane l'existence d'une relation $a_1 = a_2+k$ avec $k$ dans $\cyclique$. Nous adopterons la formulation quivalente suivante, moins pesante.

\begin{definition}
\label{def:orient:cyclique}
Un graphe $\Gr$ est dit \emph{cycliquement orient aux sommets} s'il est muni d'une action de $\cyclique$ sur $\Gr_1$ dont les orbites concident exactement avec les toiles du graphe.
\end{definition}

\begin{remarque}
Il reviendrait au mme de munir indpendament chaque toile d'une action transitive de $\cyclique$ mais les notations seraient alors moins transparentes.
\end{remarque}

\begin{definition}
Un \emph{morphisme} de graphes cycliquement orients aux sommets est un morphisme de graphes $\phi$ dont la composante $\phi_1$ est un morphisme de $\cyclique$-ensembles (\idest $\phi_1(a+1)=\phi_1(a)+1$ quel que soit $a$).
\end{definition}

\setcounter{principe}{0}

\begin{principe}
\label{princ:action:produit:libre}
La catgorie des ensembles munis de deux actions indpendantes, d'un groupe $G_1$ et d'un groupe $G_2$, o les morphismes sont les applications qui commutent simultanment aux deux actions, est \emph{isomorphe}  la catgorie des ensembles munis d'une action du produit libre $G_1 \ast G_2$.
\end{principe}

\begin{demonstration}
Cela rsulte de la proprit universelle du produit libre : en notant $\rho$ les morphismes naturels associs au produit libre et $\eta$ les morphismes de structure de ces actions de groupes, on voit qu'ils s'insrent dans un diagramme commutatif comme suit : 
\begin{align*}
\xymatrix@C=.5cm@R=.7cm
{
	& {G_1 \ast G_2}
		\ar@{..>}[dd]^{\eta_1\ast\eta_2}
	\\ {G_1}
		\ar[ur]^{\rho_1}
		\ar[dr]_{\eta_1}
	&& {G_2}
		\ar[ul]_{\rho_2}
		\ar[dl]^{\eta_2}
	\\ & {\Sym_X}
}
\end{align*}
Fin de la dmonstration.
\hfill $\Box$ \end{demonstration}

\begin{application*}
En faisant $G_1 = \cyclique$ et $G_2 = \cycliquedeux$, on voit qu'un graphe $\Gr$ cycliquement orient aux sommets est canoniquement muni d'une action du groupe $\graphique$, produit libre du groupe cyclique infini $\cyclique$ avec le groupe  deux lments $\cycliquedeux$, sur l'ensemble de ses arcs, qui soit compatible aux morphismes naturels $\cyclique \to \graphique \leftarrow \cycliquedeux$.
\end{application*} \shortpar

L'opration qui  un graphe muni d'une orientation cyclique aux sommets associe l'ensemble de ses arcs munis de cette action de groupe est fonctorielle (\emph{foncteur d'oubli}). On voit tout de suite qu'elle ne peut pas tre fidle puisque deux graphes sans artes mais comportant un nombre de sommets diffrent ne sont pas distingus  l'arrive. On est donc amen  restreindre la catgorie de \emph{dpart} :

\begin{definition}
\label{def:diagramme}
On dira qu'un graphe muni d'une orientation cyclique aux sommets est un \emph{diagramme} s'il ne comporte aucun sommet isol.
\end{definition}

\begin{variante*}
Un diagramme sera qualifi de \emph{trivalent} si ses sommets sont tous d'ordre \emph{un} ou \emph{trois}.
\end{variante*}

\begin{lemme}
\label{lemme:diag:triv}
Un diagramme est trivalent si et seulement si l'action de $\ZZ$ sur ses arcs vrifie la condition $a+3=a$ pour tout $a$.
\end{lemme}

\begin{demonstration}
C'est clair.
\end{demonstration}

\begin{theoreme}
\label{equiv:cat:diag:gens}
Le foncteur d'oubli, qui  un diagramme associe l'ensemble de ses arcs muni de l'action du groupe $\graphique$, est une quivalence de catgories.
\end{theoreme}

\begin{demonstration}
Nous allons dcrire une opration de \emph{reconstruction} qui  un ensemble muni d'une action de $\graphique$ redonne un diagramme \emph{naturellement isomorphe} au diagramme de dpart. Nous n'expliciterons pas la fonctorialit car celle-ci dcoule immdiatement de la nature des oprations utilises.

Soit $X$ un ensemble muni d'une action du groupe $\graphique$, on lui associe le graphe dont l'ensemble des arcs est $X$ et dont l'ensemble des sommets est le quotient de $X$ par l'action du groupe $\cyclique$. L'application $s$ est la projection canonique $X \to X/\cyclique$, et l'application $.\inv : X \to X$ est l'involution associe  l'lment non-trivial de $\cycliquedeux$. L'application $t : X \to X/\cyclique$ est alors dfinie par $t(a) = s(a\inv)$.

Les trois conditions de la dfinition \ref{def:graphes} sont bien vrifies.  D'autre part, vu la construction, l'action de $\ZZ$ est bien compatible  l'application $s$ au sens de la dfinition \ref{def:orient:cyclique}. La dfinition \ref{def:diagramme} est elle aussi satisfaite puisque l'application $s$ est \emph{surjective}.

Il reste maintenant  exhiber l'isomorphisme naturel $\phi$ entre le diagramme de dpart $\Gr$ et le diagramme reconstruit via cette opration :
pour la composante $\phi_1$ on choisit simplement l'identit de $\Gr_1$, 
la composante $\phi_0$ ncessite en revanche un raisonnement. Il est immdiat au vu des dfinitions \ref{def:orient:cyclique} et \ref{def:diagramme} que la paire $(\Gr_0,s)$ est un reprsentant du foncteur qui  tout ensemble $Y$ associe l'ensemble des applications $u : \Gr_1 \to Y$ pour lesquelles on ait $u(a+1)=u(a)$ quel que soit $a\in \Gr_1$. Il rsulte un isomorphisme naturel \emph{canonique} entre la paire $(\Gr_0,s)$ et la paire $(\Gr_1/\ZZ,\pi)$ o $\pi$ dsigne la projection $\Gr_1 \to \Gr_1/\ZZ$. On fait le choix de cette bijection pour $\phi_0$. La commutativit des trois diagrammes de la dfinition \ref{def:graphes:morph} est alors immdiate par naturalit. Fin de la dmonstration.
\hfill $\Box$ \end{demonstration}

Le rsultat suivant permet de raffiner l'quivalence que l'on vient d'obtenir.

\begin{theoreme}
\label{connexe:equiv:transitive}
Un diagramme est connexe si et seulement si l'action du groupe $\graphique$ sur ses arcs est transitive.
\end{theoreme}

\begin{demonstration}
On rappelle qu'un graphe est connexe s'il existe une chane reliant toute paire de sommets distincts $x$ et $y$. C'est--dire une suite finie d'arcs $a_0,...,a_n$ telle que $s(a_0) = x$, $t(a_n)=y$ et $t(a_{k-1}) = s(a_{k})$ pour $k = 1,...,n$. Comme l'action de $\ZZ$ est suppose transitive sur chaque toile, il revient au mme de supposer qu'il existe une suite d'entiers $r_1,...,r_n$ telle que $a_{k-1}\inv = a_{k} + r_k$. Comme d'autre part l'application $s$ est suppose surjective, cela quivaut bien  ce que l'action de $\graphique$ soit transitive sur les arcs. Fin de la dmonstration.
\hfill $\Box$ \end{demonstration}

\begin{corollaire*}
\label{equiv:cat:diag:con:gens:trans}
Il rsulte que la catgorie des diagrammes connexes est quivalente  celle des $\graphique$-ensembles transitifs.
\end{corollaire*}

\begin{convention*}
Dans ce qui suit, sauf mention expresse du contraire, les diagrammes seront tous supposs \emph{connexes}.
\end{convention*}

\subsection{Diagrammes Points}

\begin{definition}
Un diagramme est dit \emph{point} lorsqu'il est muni d'un \emph{arc distingu}, lequel est appel \emph{point base} du diagramme (\emph{sic}).
Les morphismes de diagrammes points seront toujours supposs respecter les points base.
\end{definition}
\shortpar
Cette dfinition est justifie par le rsultat suivant :

\begin{lemme}
Un morphisme $\phi$ entre deux diagrammes connexes $\Gr$ et $\Gr'$ est compltement dtermin par la donne d'un arc de $\Gr$ et de son image.
\end{lemme}

\shortpar
Il est clair qu'inversement, il n'existe pas toujours de morphisme d'un diagramme point $\Gr$ dans un autre $\Gr'$. On dispose d'un critre prcis :

\begin{lemme}
\label{morph:point:fix}
Il existe un morphisme point $\phi : (\Gr,a) \to (\Gr',a')$ entre deux diagrammes connexes $\Gr$ et $\Gr'$ si et seulement si les fixateurs de $a$ et de $a'$ dans $\graphique$ satisfont  la relation d'inclusion suivante {\rm:}
\begin{align*}
	\Fix_a(\graphique) \;\subseteq\; \Fix_{a'}(\graphique)
\end{align*}
\end{lemme}

\shortpar
En vertu de l'quivalence de catgories dmontre au paragraphe prcdent, les deux lemmes sont une consquence immdiate du principe gnral suivant~:

\begin{principe}
\label{fix:morph:pinc}
Etant donns deux ensembles points $(X,x)$ et $(Y,y)$ munis chacun d'une action transitive d'un mme groupe $G$, il existe une application quivariante pointe $f : (X,x) \to (Y,y)$ si et seulement si les fixateurs de $x$ et de $y$ dans $G$ satisfont  la relation d'inclusion suivante {\rm:}
\begin{align*}
	\Fix_x(G) \;\subseteq\; \Fix_{y}(G)
\end{align*}
Une telle application est ncessairement unique.
\end{principe}

\begin{remarque}
Il n'est pas ncessaire de supposer que le groupe $G$ opre transitivement sur l'ensemble d'arriv $Y$. Toutefois, si l'on ne suppose plus que l'action de $G$ sur $X$ soit transitive, il faut gnraliser la notion de pointage pour retrouver un nonc analogue.
\end{remarque}

\begin{demonstration}

Comme l'action de $G$ sur $X$ est suppose transitive, la condition que $f$ soit pointe et quivariante suffit  la dterminer compltement. En effet, pour tout lment $x'$ de $X$, on peut choisir un lment $g$ de $G$ tel que $x'=g.x$ et alors $f(x')=f(g.x) = g.f(x) = g.y$.
Montrons  prsent que la fonction $f$ est indpendante de ce choix.
Cela correspond  ce que
$g.x = g'.x$ entrane $g.y = g'.y$ quels que soient $g$ et $g'$ deux lments de $G$. Il est par ailleurs vident que cette condition revient exactement  ce que les fixateurs de $x$ et de $y$ vrifient l'inclusion de l'nonc. Fin de la dmonstration.
\hfill $\Box$
\end{demonstration}

\begin{corollaire*}
Une consquence importante de l'unicit que l'on observe dans le principe prcdent, est que si l'on a deux morphismes points comme suit,
\begin{align*}
	(\Gr,a) \to (\Gr',a')
	\quad \text{ et } \quad
	(\Gr',a') \to (\Gr, a)
\end{align*}
ce sont ncessairement deux isomorphismes rciproques l'un de l'autre.
\end{corollaire*}

\subsection{Aspects Algorithmiques}
Une consquence intressante de l'enrichissement que constitue la donne d'une orientation cyclique en chaque sommet, est qu'il est trs commode de dcider \emph{rcursivement} (c'est--dire au moyen d'un \emph{algorithme}) s'il existe un morphisme point $\phi : (\Gr,a) \to (\Gr',a')$ entre deux diagrammes connexes \emph{finis} $\Gr$ et $\Gr'$.

\begin{definition}
\label{paires:critiques}
Une \emph{paire critique} (relativement  deux diagrammes $\Gr$ et $\Gr'$) est par dfinition un couple $(a,a')$ constitu d'un arc de chacun des graphes ($a$ dans $\Gr_1$ et $a'$ dans $\Gr'_1$) pour lesquels il existe une \emph{obstruction}  la ralisation d'un morphisme point $(\Gr,a)\to(\Gr',a')$.
\end{definition}
\shortpar
Une telle obstruction prend en vertu du lemme \ref{morph:point:fix} ci-dessus, la forme d'un lment $g$ de $\graphique$ tel que $g.a = a$ et $g.a'\neq a'$.
Une faon pratique de dcider de l'existence d'un tel lment consiste par exemple,  calculer la clture dans $\Gr_1\times \Gr_1'$ de la paire $(a,a')$ par les deux rgles suivantes~:
\begin{align*}
	(a_1,a_2) \quad&\longrightarrow \quad (a_1 + 1, a_2 + 1) \\
	(a_1,a_2) \quad&\longrightarrow \quad (a_1 \inv, a_2 \inv)
\end{align*}
Il s'agit de la plus petite partie de $\Gr_1\times \Gr_1'$ qui contienne $(a,a')$ et qui soit stable vis--vis de ces deux rgles.

Comme l'ensemble $\Gr_1\times \Gr_1'$ est \emph{fini} le calcul \emph{termine} et d'aprs le lemme, la paire $(a,a')$ considre initialement est critique si et seulement si l'on rencontre deux paires $(a_1,a_2)$ et $(a_3,a_4)$ avec $a_1 = a_3$ et $a_2 \neq a_4$. En arrtant le calcul ds que cette situation se produit,
on peut mme relcher la condition de finitude portant sur le graphe $\Gr'$ ; on voit en outre que le nombre d'oprations lmentaires est du mme ordre de grandeur que le nombre d'arcs de $\Gr$.


\section{Principe de Classification}

La classification repose en premier lieu sur l'quivalence de catgories entre les diagrammes points et les $G$-ensembles transitifs points faisant l'objet du thorme \ref{equiv:cat:diag:gens}. Il convient ici d'noncer quelques gnralits supplmentaires concernant cette dernire catgorie.

\begin{principe}
\label{G:ens:pont:libre:tans}
L'ensemble sous-jacent  un groupe $G$ muni de son action par translation {\rm(} gauche{\rm)} et point par l'lment neutre est un objet \emph{initial} dans la catgorie des $G$-ensembles transitifs points Ñ C'est en ralit le cas de tout ensemble point sur lequel $G$ opre librement et transitivement.
Le point muni de l'action triviale du groupe $G$ et de son pointage tautologique est quant  lui \emph{terminal}.
\end{principe}

\begin{principe}
\label{princ:coset:action}
Soit $H$ un sous-groupe quelconque d'un groupe $G$ et soit $(X,x)$ un ensemble point sur lequel $G$ opre \emph{ gauche}, librement et transitivement {\rm(}voir par exemple le principe \ref{G:ens:pont:libre:tans} ci-dessus{\rm)}.
On dfinit une action \emph{ droite} de $H$ sur $X$ par,
\begin{align*}
	x' . h \overset{\text{def.}}{=} (gh . x)
	\qquad \text{Pour tout $x' = g . x \in X$ et tout $h \in H$.}
\end{align*}
Cette action dpend videmment du point base $x$.
On note $X/H$ l'ensemble quotient et $\pi : X \to X/H$ la projection correspondante.
Le quotient $X/H$ est alors canoniquement muni d'une action  gauche de $G$ pour laquelle la projection $\pi$ soit quivariante.
Cette action est transitive et de plus, le fixateur de $\pi(x)$ dans $G$ est exactement $H$.
\end{principe}

\begin{remarque}
En vertu du principe \ref{fix:morph:pinc}, ce qui prcde est un problme universel. L'objet construit est donc unique  isomorphisme canonique prs.
\end{remarque}

\begin{demonstration}
L'action de $H$  droite est bien dfinie puisque par hypothse la reprsentation $x' = g.x$ existe ($G$ agit transitivement) et est unique ($G$ agit librement).
Vrifions que la relation $g . \pi(x') = \pi(g.x')$ dfinit une action de $G$ sur le quotient $X/H$, ou plus prcisment, que cette relation ne dpend pas du reprsentant $x'$ de la classe $\pi(x')$. Soit donc $x'$ et $x''$ deux lments de $X$ pour lesquels $\pi(x') = \pi(x'')$, il existe donc $h\in H$ tel que $x'.h = x''$. On pose $x' = g'.x$ et $x'' = g'' . x$ avec $g'$ et $g''$ dans $G$ (\emph{uniques}, comme ci-dessus). Il rsulte que l'on a $g' h . x = g''. x$. En faisant agir l'lment $g$  gauche on obtient $gg' h.x= gg''.x$. En appliquant $\pi$ aux deux membres de l'galit on obtient que $\pi(gg'h.x) = \pi(gg''.x)$. En faisant agir l'lment $h\inv$  droite sur l'argument de $\pi$ dans le membre de gauche, on obtient finalement la relation $\pi(gg'.x) = \pi(gg''.x)$ d'o $\pi(g.x') = \pi(g.x'')$. Ceci achve de dmontrer que l'action de $G$ induite par $\pi$ est bien dfinie.
La dernire assertion est vidente.
Fin de la dmonstration.
\hfill $\Box$ \end{demonstration}

\begin{principe}
\label{principe:fix:conj}
On sait bien que dans un $G$-ensemble, les fixateurs de deux lments d'une mme orbite sont conjugus. Plus prcisment, on a :
\begin{align*}
	\Fix_x(G) = g\cdot\,\Fix_{g.x}(G)\,\cdot g\inv
\end{align*}
\end{principe}

\begin{application}
Le principe \ref{princ:coset:action} ci-dessus achve de dmontrer les deux correspondances biunivoques suivantes :
\begin{center}
\begin{tabular}{ccccc}
	\begin{tabular}{c}
	Diagrammes points \\  isomorphismes prs
	\end{tabular}
	& $\longleftrightarrow$ &
	\begin{tabular}{c}
	Sous-groupes de $\ZZ\ast\ZZ/2\ZZ$
	\end{tabular}
	\\ \\
	\begin{tabular}{c}
	Diagrammes trivalents points \\  isomorphismes prs
	\end{tabular}
	& $\longleftrightarrow$ &  
	\begin{tabular}{c}
	Sous-groupes de $\mathrm{PSL}_2(\ZZ)$
	\end{tabular}
	\\ \\
\end{tabular}
\end{center}
\end{application}

\begin{application}
Le principe \ref{principe:fix:conj} ci-dessus permet de \emph{gommer} les points base. Il en rsulte les deux quivalences suivantes.
\begin{center}
\begin{tabular}{ccccc}
	\begin{tabular}{c}
	Diagrammes \\  isomorphismes prs
	\end{tabular}
	& \hspace{.5cm} & $\longleftrightarrow$ & \hspace{.5cm}
	\begin{tabular}{c}
	Sous-groupes de $\ZZ\ast\ZZ/2\ZZ$ \\  conjugaison prs
	\end{tabular}
	\\ \\
	\begin{tabular}{c}
	Diagrammes trivalents \\  isomorphismes prs
	\end{tabular}
	& \hspace{.5cm} & $\longleftrightarrow$ & \hspace{.5cm} 
	\begin{tabular}{c}
	Sous-groupes de $\mathrm{PSL}_2(\ZZ)$ \\  conjugaison prs
	\end{tabular}
	\\
\end{tabular}
\end{center}
Autrement dit, deux diagrammes points correspondent  deux sous-groupes conjugus si et seulement si les diagrammes obtenus en gommant les points base sont isomorphes.
\end{application}

\begin{corollaire*}
Un diagramme correspond donc  un sous-groupe distingu si et seulement si son groupe d'automorphismes opre transitivement sur ses arcs.
\end{corollaire*}

\begin{exemples}

Les tables \ref{tab:diag:3:taille:inf:5}  \ref{tab:diag:3:taille:9} ci-contre, donnent une liste exhaustive des diagrammes trivalents connexes dont la taille est infrieure  neuf. Afin de permettre une meilleure lisibilit de leurs arcs, les diagrammes sont reprsents sous la forme de leur subdivision barycentrique enrichie (\confer section \ref{sec:subdivision:barycentrique}). Enfin, dans cette reprsentation, l'orientation cyclique aux sommets est choisie pour concider avec l'orientation trigonomtrique du plan du dessin.

\begin{table}
\begin{center}
\includegraphics{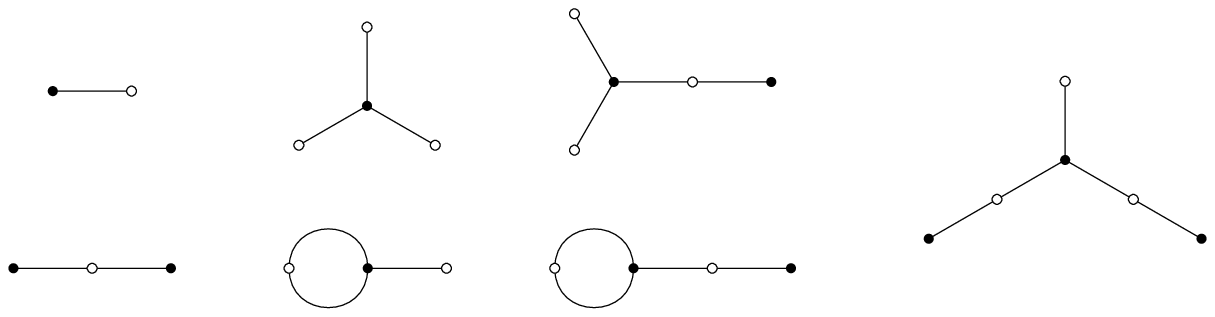}
\caption{Diagrammes trivalents de taille infrieure  cinq.}
\label{tab:diag:3:taille:inf:5}
\end{center}
\end{table}

\begin{table}
\begin{center}
\includegraphics{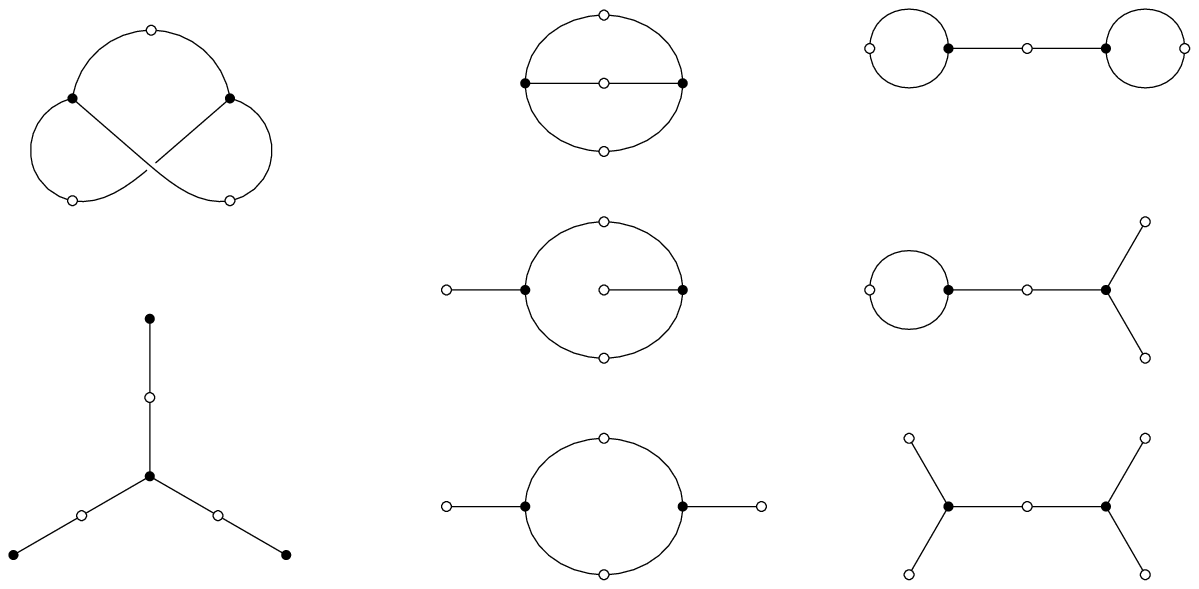}
\caption{Diagrammes trivalents de taille six.}
\label{tab:diag:3:taille:6}
\end{center}
\end{table}

\begin{table}
\begin{center}
\includegraphics{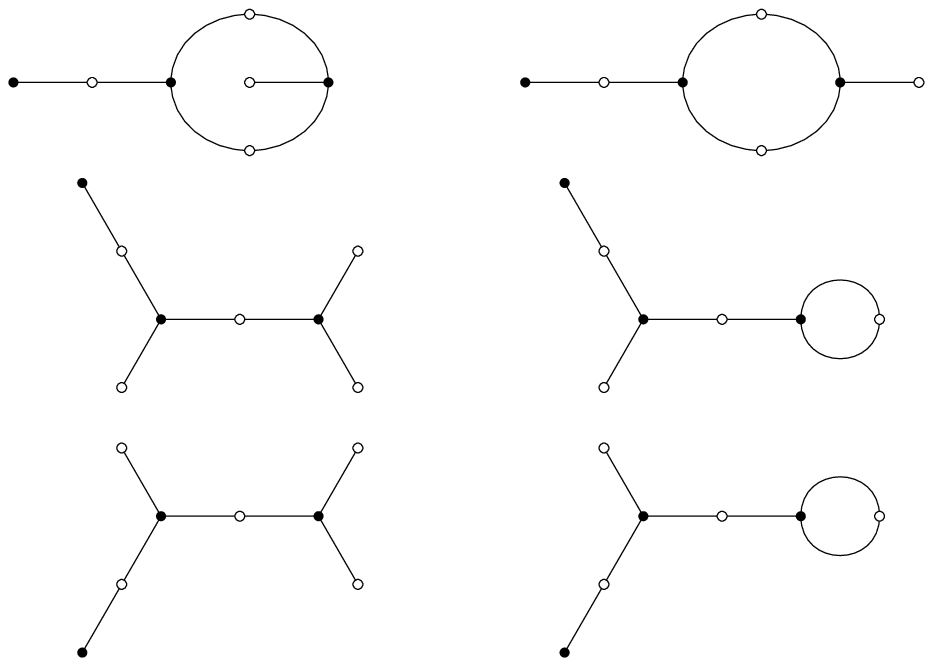}
\caption{Diagrammes trivalents de taille sept.}
\label{tab:diag:3:taille:7}
\end{center}
\end{table}

\begin{table}
\begin{center}
\includegraphics{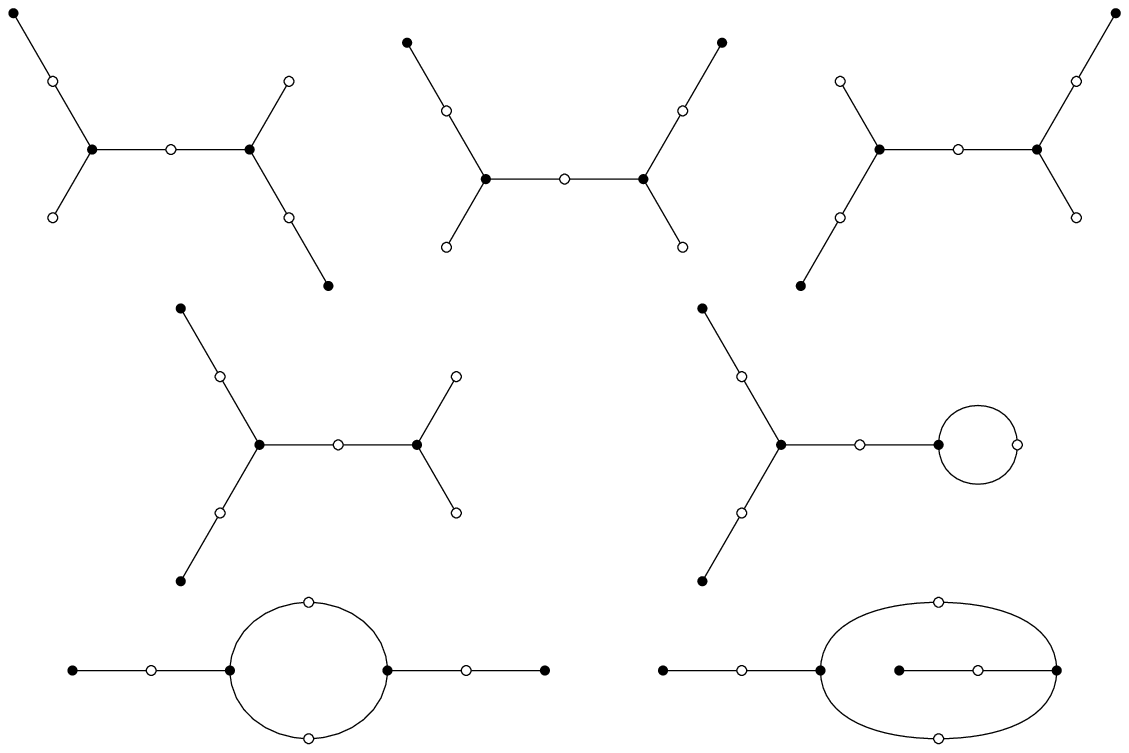}
\caption{Diagrammes trivalents de taille huit.}
\label{tab:diag:3:taille:8}
\end{center}
\end{table}

\begin{table}
\begin{center}
\includegraphics{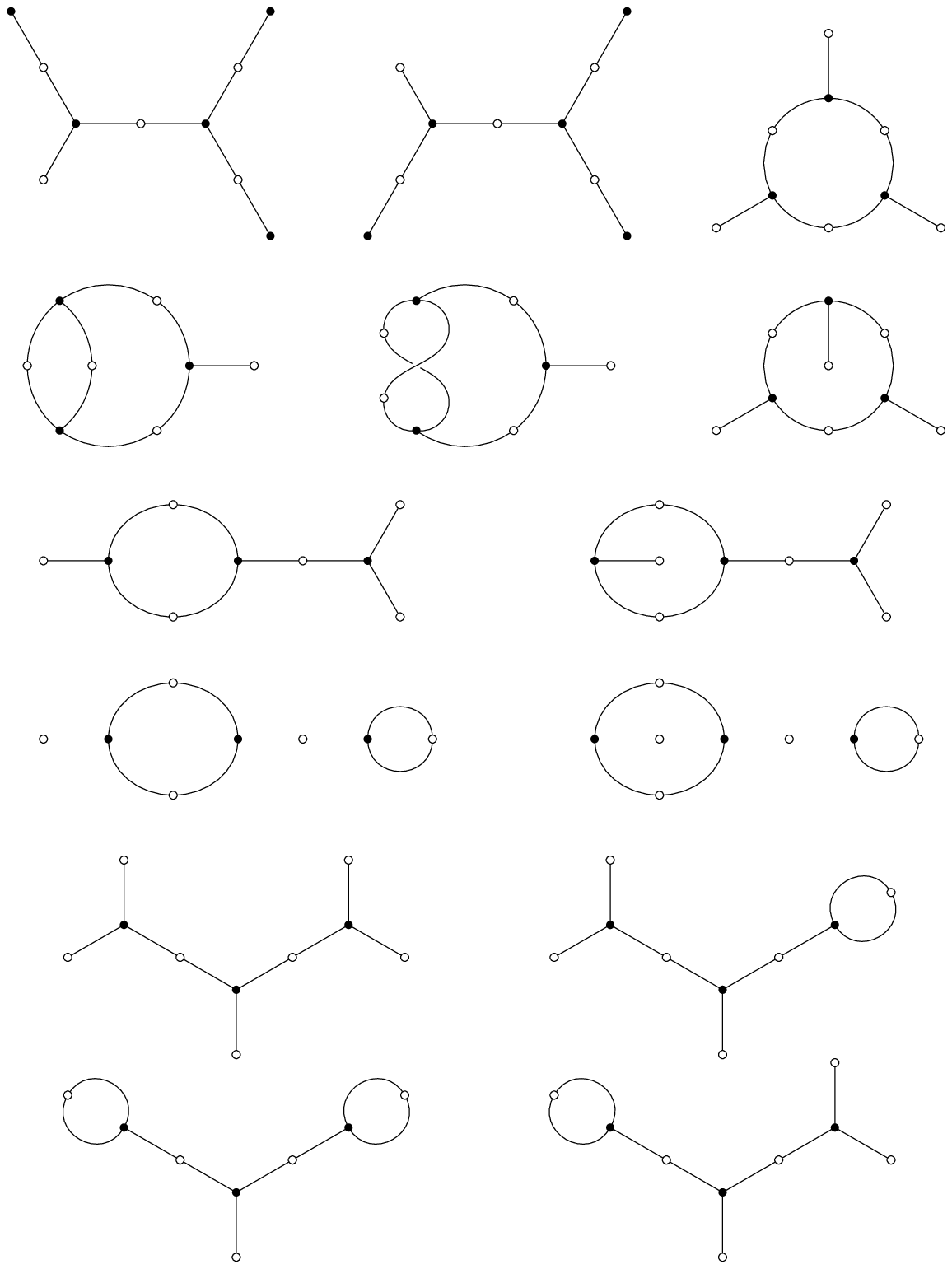}
\caption{Diagrammes trivalents de taille neuf.}
\label{tab:diag:3:taille:9}
\end{center}
\end{table}

\end{exemples}

\section{Dnombrement}

Un des objectifs de cet article est  prsent, de donner des formules gnrales, sous forme de sries gnratrices, permettant de compter le nombre de diagrammes trivalents connexes ainsi que leurs variantes pointes. Pour obtenir ces sries nous nous sommes appuys sur quelques rsultats de \emph{thorie des espces}. Par commodit pour le lecteur, nous avons cru bon de prsenter dans les sections suivantes un rapide survol des quelques rudiments qui nous sont utiles.

Le lecteur souhaitant approfondir le sujet pourra se rfrer  l'article fondateur \cite{joyal81}, les articles suivants sont consacrs  certains perfectionnements et gnralisations~: \cite{labelle81}, \cite{yeh86} et \cite{yeh87}. Enfin, l'ouvrage \cite{bergeron94} reprend remarquablement l'ensemble des dfinitions, notations, et rsultats de cette thorie, en un expos  la fois pdagogique et synthtique.

\subsection{Espces Combinatoires, Gnralits}

La thorie des espces dbute par une reformulation claire de la thorie de Plya \cite{polya37, polyaread87}. Il s'agit d'encoder dans une srie dite d'index cyclique  une infinit de variables commutatives $x_1$, ...,  $x_k$, ... la combinatoire d'une famille de structures.

\begin{definition}
\label{espece:def}
Une \emph{espce combinatoire} (au sens de A. Joyal \cite{joyal81}) est un foncteur $F$ de la catgorie des ensembles finis avec bijections, vers celle des ensembles finis avec applications quelconques.
\end{definition}
\shortpar

Etant donns un ensemble fini $E$ et une espce $F$, les lments de l'ensemble $F(E)$ sont les \emph{$F$-structures}, ou \emph{structures d'espce $F$}, construites sur $E$.  On dit aussi que l'ensemble $E$ \emph{tiquette} les structures de l'ensemble $F(E)$. Si $\varrho$ est une bijection entre deux ensembles finis $E_1$ et $E_2$, l'application induite $F(\varrho)$, note $\varrho_*$, est le \emph{transport de structure} le long de $\varrho$, elle associe  toute $F$-structure sur $E_1$ une $F$-structure sur $E_2$, elle est de plus \emph{bijective} par fonctorialit. On considre cette application comme un  \emph{rtiquetage}.

\begin{definition}
\label{espece:morph:def}
Un \emph{morphisme} $\phi$ entre deux espces combinatoires $F$ et $G$ est une \emph{opration} qui  tout ensemble fini $E$ associe une application $\phi_E : F(E) \to G(E)$ commutant aux applications de transfert de structure au sens o, quel que soit l'application bijective $\varrho$ entre deux ensembles finis $E_1$ et $E_2$, le diagramme suivant soit commutatif~:
\begin{align*}
\xymatrix@C=1.5cm@R=1.5cm
{
	F(E_1)
		\ar[r]^{\phi_{E_1}}
		\ar[d]_{\varrho_*}	&
	G(E_1)
		\ar[d]^{\varrho_*}
	\\
	F(E_2)
		\ar[r]_{\phi_{E_2}}	&
	G(E_2)
}
\end{align*}
\end{definition}

On rsume la situation en disant que $\phi$ est une application \emph{naturelle} entre l'ensemble des $F$-structures construites au dessus d'un ensemble donn et celui des $G$-structures construites au-dessus du mme ensemble.

Le simple fait d'avoir dgag et mis en lumire cette condition de naturalit, doit tre considr comme un progrs important et un apport significatif de la thorie des espces. Un tel morphisme, prserve non seulement le dcompte des structures tiquetes mais aussi celui des classes de structures obtenues par l'opration de \emph{rtiquetage} \idest l'opration du groupe symtrique, induite par fonctorialit  partir de son action sur l'ensemble des tiquettes.

\begin{exemples}
Nous donnons ici quelques exemples qui seront rutiliss ensuite pour traiter les problmes numratifs viss~:
\begin{enumerate}[\qquad 1.\quad]
	\item L'\emph{espce des ensembles} note $\mathrm{Ens}$ est le foncteur qui,  tout ensemble  fini $E$, associe le singleton $\mathrm{Ens}\,(E)=\{\,E\,\}$ et  toute application bijective $\varrho : E_1 \to E_2$ associe l'application canonique $\mathrm{Ens}\,(E_1) \to \mathrm{Ens}\,(E_2)$.
	\item L'\emph{espce des permutations} note $S$ est le foncteur qui,  tout ensemble fini $E$, associe l'ensemble $S(E)$ 
	des permutations $\sigma$ de $E$ et qui  toute bijection $\varrho : E_1 \to E_2$ associe la conjugaison par $\varrho$ dfinie par $\varrho_*(\sigma) =\varrho\,\sigma\,\varrho\inv$.
	\item L'\emph{espce des permutations d'ordre $n$} note $S_n$ qui,  tout ensemble fini $E$, associe le sous-ensemble $S_n(E)$ de $S(E)$ constitu des permutation $\sigma$ vrifiant $\sigma^n = \id$.
	\item L'\emph{espce des cycles} note $C$ qui  tout ensemble fini $E$ associe le sous-ensemble $C(E)$ de $S(E)$ constitu des permutations cycliques.
	\item On dfinit de mme l'\emph{espce des cycles de longueur $n$} note $C_n$ qui,  tout ensemble fini $E$, associe l'ensemble $C_n(E) = C(E)$ si $E$ est de taille $n$ et l'ensemble vide sinon.
\end{enumerate}

\end{exemples}

\subsubsection{Structures tiquetes \textit{vs.} Stuctures Non-tiquetes}

On distingue essentiellement deux problmatiques de dnombrement attaches  une espce combinatoire donne $F$.
La premire, consiste  exprimer le nombre de structures \emph{tiquetes} d'espce $F$ construites sur un ensemble d'tiquettes donn $E$.
N'importe quelle bijection entre les deux ensembles $F(E)$ et $G(E)$ prserve ce dnombrement.

La deuxime problmatique, beaucoup plus dlicate, consiste  compter les structures \emph{non-tiquetes}. Il s'agit des orbites de l'action du groupe symtrique, induite par fonctorialit  partir de son action sur l'ensemble des tiquettes. Deux structures tiquetes tant considres comme identiques du point de vue non-tiquet ds que l'on peut obtenir l'une  partir de l'autre par une opration de \emph{rtiquetage}. Pour qu'une bijection entre les deux ensemble $F(E)$ et $G(E)$ prserve ce dnombrement il suffit qu'elle commute aux rtiquetages.

\subsubsection{Sries associes}
\label{sec:series:associees}

A une espce combinatoire $F$, on associe trois sries gnratrices $F(t)$, $\tilde{F}(t)$ et $\ZSeries_F(x_1,x_2,\dots,x_k,\dots)$ en l'indtermine commutative $t$ (de poids \emph{un}) pour les deux premires et en l'infinit d'indtermines commutatives $x_1, \dots , x_k, \dots$ (de poids $1, \dots, k, \dots$ respectivement) pour la dernire. Elles correspondent chacune  une problmatique de dnombrement prcise.

\begin{enumerate}[\qquad 1.\quad]
\item La \emph{srie gnratrice exponentielle} ou \emph{srie de Hurwitz}
\begin{equation*}
	F(t)=\sum_{n\ge 0}\frac{a_n}{n!}\, t^n\qquad \text{avec } a_n \overset{\text{def.}}{=} \abs{F(\{\,1,\dots,n\,\})}
\end{equation*}
qui compte \emph{les structures tiquetes}.
\item La \emph{srie gnratrice des types d'isomorphisme}
\begin{equation*}
	\tilde{F}(t) = \sum_{n\ge 0} \tilde{a}_n\,t^n \qquad \text{avec } \tilde{a}_n \overset{\text{def.}}{=} \abs{F(\{\,1,\dots,n\,\})/\mathfrak{S}_n}
\end{equation*}
qui compte \emph{les classes d'isomorphisme des structures}.
\item Enfin, la \emph{srie indicatrice des cycles}, ou \emph{srie de Joyal-Plya} \cite{polya37, polyaread87,joyal81},
\begin{align*}
	\ZSeries_F &= \sum_{n\ge 0}\, \frac{1}{n!}\, \sum_{\sigma \in \mathfrak{S}_n} 
			\,a_{\sigma_1,\dots,\sigma_n}\,x_1^{\sigma_1}\cdots x_n^{\sigma_n} \\
	&= \sum_{n\ge 0}\, \sum_{k_1+2k_2+\dots+nk_n=n}
			\frac{a_{k_1,\dots,k_n}}{1^{k_1}\,k_1!\cdots n^{k_n}\,k_n!}\,x_1^{k_1}\cdots x_n^{k_n}
\end{align*}
o l'on dsigne par $\sigma_1,\dots,\sigma_n$ les nombres de cycles de longueur $1$  $n$ de la permutation $\sigma$ et o $a_{\sigma_1,\dots,\sigma_n}$ dsigne le nombre de points fixes de la permutation $\sigma_* = F(\sigma)$ induite par $F$ sur l'ensemble $F(\{\,1,\dots,n\,\})$ des structures d'espce $F$ construites au-dessus de l'ensemble $X = \{\,1,\dots,n\,\}$.
La notation prcdente est justifie puisque ce nombre ne dpend que du type cyclique de $\sigma$ (\idest les entiers $\sigma_1, \dots,\sigma_n$). On peut donc utiliser la notation $a_{k_1,\dots,k_n}$, o $k_1, \dots, k_n$ sont des entiers vrifiant $k_1+2k_2+\dots+nk_n=n$.
\end{enumerate}

Les trois sries prcdentes sont des \emph{invariants} de l'espce combinatoire $F$ pour la notion d'isomorphisme d'espces combinatoires. Mais si l'on omettait la condition de naturalit prsente dans la dfinition \ref{espece:morph:def}, seule la premire srie serait conserve par isomorphisme. C'est l, la lacune principale de la combinatoire bijective qui empchait, avant la thorie des espces, de traiter des exemples significatifs de dnombrements non-tiquets. Les seuls dnombrements non-tiquets s'obtenaient tous jusque-l,  partir d'un dnombrement tiquet, en rsolvant les symtries par une astuce.

Le dnombrement de diagrammes points que l'on donne la section \ref{denombrement:diagrammes:pointes} est un exemple d'application d'une telle mthode~; il s'obtient sans recourir au sries indicatrices. Par contraste, il semble impossible d'obtenir le dnombrement de diagrammes non-points  partir d'un dnombrement de structures tiquetes. Le recours aux sries indicatrices est donc indispensable dans ce cas.

\theoremstyle{plain}
\newtheorem*{lemmedecondensation}{Lemme de condensation}

\begin{lemmedecondensation}
\label{ZF:to:Ft:and:Ftildet}
La srie indicatrice des cycles $\ZSeries_F$ d'une espce $F$ raffine les deux sries gnratrices $F(t)$ et $\tilde{F}(t)$ au sens o l'on a
$ F(t) = \ZSeries_F(t,0,0,\dots)$ et $\tilde{F}(t) = \ZSeries_F(t, t^2, t^3, \dots )$.
\end{lemmedecondensation}

\begin{remarque}
Ces deux changements de variable sont compatibles avec la graduation pose prcdemment.
\end{remarque}

\begin{demonstration} La premire formule rsulte immdiatement de l'interprtation combinatoire des coefficients. La deuxime est quant  elle  consquence immdiate du lemme de Burnside.
\hfill $\Box$ \end{demonstration}

\begin{exemples}
Nous donnons ci dessous, les formules pour la srie d'index cyclique de l'espce des ensembles et de celle des permutations. Le cas des espces  $C$, $C_n$ et $S_n$ des cycles, des cycles de longueur $n$ et des permutations d'ordre $n$ respectivement, est l'objet de la section \ref{denombrement:permutation}.
\begin{enumerate}[\qquad1.\quad]
	\item Pour l'espce $\mathrm{Ens}$ des ensembles, on a trivialement $a_n=1$, $\tilde{a}_n = 1$ et $a_{k_1,\dots,k_n}=1$. D'o il rsulte aussitt,
\begin{align*}
	\mathrm{Ens}(t) &= \sum_{n \ge 0} \, \frac{1}{n!} \,t^n = \exp (t) \\
	\widetilde{\mathrm{Ens}} (t) &= \sum_{n \ge 0} \, t^n = \frac{1}{1-t} \\
	\ZSeries_{\mathrm{Ens}} (x_1,x_2,\dots) &= \sum_{n \ge 0} \sum_{k_1+2k_2+\dots+nk_n=n} \frac{1}{1^{k_1}\,k_1!\cdots n^{k_n}\,k_n!}\,x_1^{k_1}\cdots x_n^{k_n} \\
	&= \prod_{k\ge 1} \biggl(\,\sum_{n\ge 0}\, \frac{1}{k^n \,n!} \, x_k^n \biggr)\\
	&= \prod_{k\ge 1} \exp\left(\,\frac{x_k}{k}\right)
\end{align*}
\item Pour l'espce $S$ des permutations, on a trivialement $a_n = n!$, $\tilde{a}_n = p_n$, o $p_n$ est le nombre de partitions de l'entier $n$, et $a_{k_1,\dots,k_n}= 1^{k_1}\,k_1!\cdots n^{k_n}\,k_n!$. D'o il rsulte aussitt,
\begin{align*}
	S(t) &= \sum_{n \ge 0} \, t^n = \frac{1}{1-t} \\
	\tilde{S} (t) &= \sum_{n \ge 0}\, p_n\,t^n = \prod_{k\ge 1}\,\frac{1}{1-t^k} \quad \text{d'aprs Euler,}\\
	\ZSeries_{S} (x_1,x_2,\dots) &= \sum_{n \ge 0} \sum_{k_1+2k_2+\dots+nk_n=n}
	\negthickspace \negthickspace x_1^{k_1}\cdots x_n^{k_n} \\
	&= \prod_{k\ge 1} \biggl(\,\sum_{n\ge 0} \, x_k^n \biggr)\\
	&= \prod_{k\ge 1} \,\frac{1}{1-x_k}
\end{align*}
\end{enumerate}
\end{exemples}

Nous donnons dans la section suivante, les formules pour la srie d'index cyclique des espces $C$, $C_n$ et $S_n$ des cycles, des cycles de longueur $n$ et des permutations d'ordre $n$ respectivement.

\subsection{Exponentiation et Logarithme d'Espces Combinatoires}
Les formules concernant la srie d'index cyclique des trois espces $C$, $C_n$ et $S_n$ des \emph{cycles}, des \emph{cycles de longueur} $n$ et des \emph{permutations d'ordre} $n$ respectivement s'obtiennent par un jeu faisant intervenir les deux oprations d'\emph{exponentiation} et de \emph{logarithme} d'une espce combinatoire.

On suppose donne une espce combinatoire $F$ pour laquelle on ait $F(\emptyset)=\emptyset$.
On note $\mathrm{Ens}(F)$ ou encore $F^*$, l'espce obtenue  partir de l'espce $F$ par \emph{composition} plthystique avec le foncteur $\mathrm{Ens}$ de l'espce des ensembles.
Il s'agit de l'\emph{exponentielle} de l'espce $F$. Les $\mathrm{Ens}(F)$-structures sont par dfinition des \emph{assembles} \cite{joyal81} de $F$-structures. De faon plus prcise,  le foncteur $\mathrm{Ens}(F)$ associe  tout ensemble fini $E$, l'ensemble dont les lments sont les uplets $(n\,;P_1,\dots,P_n\,;A_1,\dots,A_n)$ constitus d'un entier $n$, d'une partition de $E$ comportant $n$ parties $P_1,\dots,P_n$ et sur chacune de ces parties, la donne $A_k$ d'une structure d'espce $F$ construite au dessus d'elle (\idest on impose que $A_k$ appartienne  l'ensemble $F(P_k)$ quel que soit $k$).

Pour les sries de Hurwitz associes  deux espces $F$ et $G$ relies par la relation fonctorielle $G \simeq \mathrm{Ens}(F)$, on a les relations suivantes, lesquelles justifient la terminologie.
\begin{align*}
	G(t) = \exp(F(t))
	\quad \text{ et } \quad
	F(t) = \log(G(t))
\end{align*}
Les sries gnratrices des types d'isomorphisme satisfont aux relations plus compliques suivantes, dues  Harary et Palmer \cite{hararypalmer73}.
\begin{align*}
	\tilde{G}(t) = \exp\biggl(\, \sum_{n\ge1} \frac{1}{n} \, \tilde{F}(t^n)\biggr)
	\quad \text{ et } \quad
	\tilde{F}(t) = \sum_{n\ge 1} \frac{\mu(n)}{n}\,\log\left(\tilde{G}(t^n)\right)
\end{align*}
o $\mu$ dsigne la fonction de Mbius.

Ces formules admettent la gnralisation suivante, due  Joyal \cite{joyal81}, au cas des sries d'index cyclique associes aux espces $F$ et $G$.
\begin{align*}
	\ZSeries_G = \exp\biggl(\, \sum_{k\ge 1} \frac{1}{k}\,\ZSeries_{F,k}\biggl)
	\quad\text{ et }\quad
	\ZSeries_F = \sum_{k\ge 1} \frac{\mu(k)}{k}\,\log\left(\ZSeries_{G,k}\right)
\end{align*}
o $\ZSeries_{F,k}$ et $\ZSeries_{G,k}$ dsignent les sries obtenues  partir des sries $\ZSeries_{F}$ et $\ZSeries_{G}$ respectivement, en effectuant le changement de variable consistant  remplacer la variable $x_n$ par la variable $x_{kn}$, et cela simultanment pour tout $n \ge 1$. Autrement dit, si de manire gnrale, on note $\ZSeries_{F}(x_1,x_2,x_3,\dots)$ la srie indicatrice d'une espce $F$, quel que soit l'entier $k\ge 1$, on a~:
\begin{align*}
	\ZSeries_{F,k}(x_1,x_2,x_3,\dots) &= \ZSeries_{F}(x_{k},x_{2k},x_{3k},\dots)
\end{align*}

\subsection{Produit Cartsien de deux Espces, Somme Directe}

Le produit cartsien de deux espces, s'interprte comme une \emph{superposition de structures}, contrairement au produit habituel qui lui, s'interprte comme une \emph{juxtaposition de structures}. La somme directe s'interprte elle, comme une \emph{disjonction de structures}.
De faon plus prcise :

\begin{definition}
Le \emph{produit cartsien} de deux espces combinatoires $F$ et $G$, est l'espce note $F\times G$,  dfinie pour tout ensemble fini $E$ et toute application bijective $f$ entre deux ensembles finis $E_1$ et $E_2$, par :
\begin{align*}
	(F\times G) (E) &= F(E) \times G(E)	\\
	(F \times G) (f) \, &= F(f) \, \times G(f)
\end{align*}
\end{definition}

\begin{definition} On dfinit de mme, la \emph{somme directe} de deux espces combinatoires $F$ et $G$, note $F + G$, qui  tout ensemble fini $E$ et toute application bijective $f$ entre deux ensembles finis $E_1$ et $E_2$, associe respectivement~:
\begin{align*}
	(F + G) (E) &= F(E) \sqcup G(E)\\
	(F + G) (f) \, &= F(f) \, \sqcup G(f)
\end{align*}
\end{definition}

De faon plus prcise, si l'on note :
\begin{align*}
	E_1' &= F(E_1) &
	E_2' &= F(E_2) &
	f' &= F(f) \\
	E_1'' &= G(E_1) &
	E_2'' &= G(E_2) &
	f'' &= G(f)
\end{align*}
On a deux diagrammes commutatifs comme suit :
\begin{gather*}
\xymatrix@C=1.3cm@R=1.3cm
{
	E_1'
		\ar[d]_{f'}
		&
	E_1' \times E_1''
		\ar[l]_-{\pi_1'}
		\ar[r]^-{\pi_1''}
		\ar@{..>}[d]^{f'\times f''}
		&
	E_1''
		\ar[d]^{f''}	\\
	E_2' &
	E_2' \times E_2''
		\ar[l]^-{\pi_2'}
		\ar[r]_-{\pi_2''}
		&
	E_2''
}
\qquad
\xymatrix@C=1.3cm@R=1.3cm
{
	E_1'
		\ar[d]_{f'}
		\ar[r]^-{\rho_1'}
		&
	E_1' \sqcup E_1''
		\ar@{..>}[d]^{f'\sqcup f''}
		&
	E_1''
		\ar[l]_-{\rho_1''}
		\ar[d]^{f''}	\\
	E_2' 
		\ar[r]_-{\rho_2'}
		&
	E_2' \sqcup E_2''
		&
	E_2''
		\ar[l]^-{\rho_2''}
}
\end{gather*}
o l'on note $\pi$ les projections naturelles associes  chacun des deux produits cartsiens et $\rho$ les injections naturelles associes  chacune des deux runions disjointes. De la sorte, la flche verticale du milieu de chacun des deux diagrammes est ncessairement unique et le produit cartsien de deux espces $F$ et $G$ est le \emph{produit direct}, au sens des catgories, des deux foncteurs correspondants. Idem pour la somme directe, qui s'interprte comme une \emph{somme directe} de foncteurs ou \emph{coproduit}.

\begin{remarque}
Le produit cartsien $F \times G$ de deux espces combinatoires $F$ et $G$ ne correspond pas au produit habituel, not $F\, \cdot \, G$ et il en est de mme pour les sries associes. Parmi ces deux notions, nous n'aurons toutefois  considrer que le cas du produit cartsien et non celui du produit usuel.
\end{remarque}

\subsubsection{Dnombrement tiquet}

Soient donnes deux espces combinatoires $F$ et $G$, 
dont les sries de Hurwitz sont les suivantes,
\begin{align*}
F(t) = \sum_{n \ge 0} \frac{a_n}{n!}\,t^n
\quad\text{ et } \quad
G(t) = \sum_{n \ge 0} \frac{b_n}{n!}\,t^n
\end{align*}
Il est clair que la srie de Hurwitz du produit cartsien $F \times G$ de ces deux espces est le \emph{produit de Hadamard} des deux sries prcdentes~:
\begin{align*}
F(t) \odot G(t) \overset{\text{def.}}{=} \sum_{n \ge 0} \frac{a_n\,b_n}{n!}\,t^n
\end{align*}

Il n'existe malheureusement pas de formule de ce genre pour les sries gnratrices des types d'isomorphisme. La difficult provient du fait que le foncteur de \emph{passage au quotient} (qui  un ensemble sur lequel opre un groupe fait correspondre l'ensemble des orbites de cette action de groupe) soufre d'un dfaut d'exactitude \emph{ gauche}~: il ne commute pas au produit cartsien et n'a pas d'adjoint formel  droite.

Il commute toutefois  la runion disjointe de sorte que la srie gnratrice des types d'isomorphisme associe  la somme directe de deux espces combinatoires est la somme usuelle des sries correspondantes.

\subsubsection{Dnombrement non-tiquet}

Pour traiter de la question du produit cartsien de deux espces combinatoires (correspondant  la \emph{superposition} de structures) il convient, dans le cadre du dnombrement non-tiquet, d'introduire la srie de Joyal-Plya ou \emph{srie indicatrice des cycles} de ces deux espces.
Ce type de srie compte les \emph{points fixes} de l'action de rtiquetage, et non ses \emph{orbites}. Or le foncteur de \emph{passage aux points fixes} d'une action de groupe est  la fois \emph{exact  gauche} et \emph{exact  droite}, il commute de fait, aux produits et aux sommes et admet un adjoint formel  gauche et  droite. 
Cela permet de retrouver un principe de comptage. Il rsulte en particulier, le lemme suivant qui tablit des formules qui combines au lemme de condensation permettent de traiter cette question.

\begin{lemme}
\label{lemme:produit:Hadamard}
Si l'on note par $\ZSeries_F$ et $\ZSeries_G$ les sries \emph{indicatrices des cycles} \cite{polya37, polyaread87,joyal81} de deux espces combinatoires $F$ et $G$~{\rm:}
\begin{align*}
	\ZSeries_F &= \sum_{n\ge 0}\, \sum_{k_1+2k_2+\dots+nk_n=n}
			\frac{a_{k_1,\dots,k_n}}{1^{k_1}\,k_1!\cdots n^{k_n}\,k_n!}\,x_1^{k_1}\cdots x_n^{k_n}
			\text{ et } \\
	\ZSeries_G &= \sum_{n\ge 0}\, \sum_{k_1+2k_2+\dots+nk_n=n}
			\frac{b_{k_1,\dots,k_n}}{1^{k_1}\,k_1!\cdots n^{k_n}\,k_n!}\,x_1^{k_1}\cdots x_n^{k_n}
\end{align*}
La srie indicatrice des cycles du produit cartsien $F \times G$ est simplement le \emph{produit de Hadamard} des deux sries prcdentes~{\rm :}
\begin{align*}
	\ZSeries_F \odot \ZSeries_G \overset{\text{def.}}{=} \sum_{n\ge 0}\, \sum_{k_1+2k_2+\dots+nk_n=n}
			\frac{a_{k_1,\dots,k_n}\, b_{k_1,\dots,k_n}}{1^{k_1}\,k_1!\cdots n^{k_n}\,k_n!}\,x_1^{k_1}\cdots x_n^{k_n}
\end{align*}
La srie indicatrice des cycles de la somme directe $F + G$ est la somme usuelle $\ZSeries_F+\ZSeries_G$ {\rm(}terme  terme{\rm)} des deux sries $\ZSeries_F$ et $\ZSeries_G$.
\end{lemme}

\begin{application}
\label{espece:diag:connexes}

Considrons les espces $D$ et $D^*$, des \emph{diagrammes} au sens de la dfinition \ref{def:diagramme}, respectivement connexes et non-ncessairement connexes. On rappelle qu'il s'agit de graphes au sens des dfinitions \ref{def:graphes} et \ref{def:graphes:morph} munis d'une orientation cyclique aux sommets au sens de la dfinition \ref{def:orient:cyclique} et dpourvus de sommets isols.
En vertu du thorme \ref{equiv:cat:diag:gens} on sait que l'espce $D^*$ est isomorphe  celle des ensembles munis d'une action du groupe $\graphique$. On sait aussi, ce qui est quivalent par application du principe \ref{princ:action:produit:libre}, qu'elle est isomorphe  l'espce des ensembles munis chacun d'une action indpendante des deux groupes $\cyclique$ et $\cycliquedeux$.
De l, on tire les deux isomorphismes naturels suivants :
\begin{align*}
	\mathrm{Ens}(D) &\simeq D^* \simeq S_2 \times S
\end{align*}
o l'on dsigne par $S_2$ l'espce des \emph{involutions} et par $S$ celle des \emph{permutations} quelconques et o le premier isomorphisme traduit l'existence et l'unicit de la dcomposition d'un diagramme en somme directe de ses composantes connexes.
\end{application}

\begin{application}
\label{espece:diag:trivalents:connexes}

Si l'on se restreint aux espces $D_3$ et $D_3^*$ des \emph{diagrammes trivalents}, respectivement connexes et non-ncessairement connexes, on a comme prcdemment, conformment au lemme \ref{lemme:diag:triv}, les deux isomorphismes naturels suivants~:
\begin{align*}
	\mathrm{Ens}(D_3) &\simeq D_3^* \simeq S_2 \times S_3
\end{align*}
o l'on dsigne comme ci-dessus, par $S_2$ l'espce des involutions et par $S_3$ celle des permutations d'ordre trois.
\end{application}

\subsection{Dnombrement de Diagrammes Points}
\label{denombrement:diagrammes:pointes}

Nous considrons dans cette section les deux espces $D_3^\bullet$ et $D_3$, des diagrammes trivalents connexes, respectivement \emph{points} et \emph{non-points}. Nous donnerons ci-dessous des exemples permettant de bien saisir la distinction. Nous considrons les sries formelles suivantes associes  ces deux espces.
\begin{center}
\begin{tabular}{lc|ccccc}
	&&&& \\
	Diagrammes & \quad & \quad & Points & \qquad& Non-points \\ &&& \\
	\hline &&&& \\
	Etiquets & & &
	$ \displaystyle
D_3^\bullet (t) = \sum_{n\ge 0} \, \dfrac{a_n^\bullet}{n!}\,t^n
$ & &
	$\displaystyle D_3 (t) = \sum_{n\ge 0} \, \dfrac{a_n}{n!}\,t^n$ \\ &&&&\\
	Non-Etiquets & & &
	$\displaystyle \tilde{D}_3^\bullet (t) = \sum_{n\ge 0} \, \tilde{a}_n^\bullet\,t^n$ & &
	$\displaystyle \tilde{D}_3 (t)= \sum_{n\ge 0} \, \tilde{a}_n\,t^n$
	\\ &&&&
\end{tabular}
\end{center}
On a trivialement la relation suivante, qui exprime que le choix d'un point base sur les structures tiquetes revient  appliquer l'oprateur d'Euler  la sries gnratrices de Hurwitz.
\begin{align*}
	D_3^\bullet (t) = t\,\frac{d}{d\,t}\,D_3(t)
\end{align*}
Comme d'autre part, les diagrammes trivalents connexes points sont des structures \emph{rigides}, au sens o elle n'admettent pas d'automorphisme, on a
$\tilde{D}_3^\bullet (t) = D_3^\bullet (t)$. De cela il rsulte,
\begin{align*}
	a_n^\bullet &= n\, a_n \quad \text{ et } \quad
	\tilde{a}_n^\bullet = \frac{a_n^\bullet}{n!} \quad \text{ puis } \quad
	\tilde{a}_n^\bullet = \frac{a_n}{(n-1)!}
\end{align*}
On s'est donc ramens, pour compter le nombre $\tilde{a}_n^\bullet$ de diagrammes trivalents connexes points et non-tiquets,  compter le nombre $a_n$ de diagrammes trivalents connexes tiquets.

De l'isomorphisme $\mathrm{Ens}(D_3) \simeq S_2 \times S_3$ on tire la relation $D_3(t) = \log\left( S_2(t) \odot S_3(t)\right)$. On sait d'autre part que les sries $S_2(t)$ et $S_3(t)$ admettent les expressions simples suivantes
\begin{align*}
	S_2(t) = \exp\left(t+\frac{t^2}{2}\right)
	\quad \text{ et } \quad 
	S_3(t) = \exp\left(t+\frac{t^3}{3}\right)
\end{align*}
On obtient alors l'expression suivante permettant de calculer les coefficients de la srie gnratrice $\tilde{D}_3^\bullet(t)$~:
\begin{align*}
	D_3(t) = \log\Biggl(\,
		\sum_{n\ge0} \,n! \,t^n\sum_{\substack{n_1 + 2n_2 = n\\ n_3+3n_4=n}}
		\frac{1}{n_1!\,n_2!\,n_3!\,n_4!\,2^{n_2}\,3^{n_4}}
	\Biggr)
\end{align*}

Les coefficients de la srie sous le $\log$, laquelle n'est autre que $D_3^*(t)$, peuvent tre calculs efficacement en utilisant  l'quation de rcurrence linaire donne dans la table \ref{tab:rec:diag:3:non:connexes}. On aboutit  finalement au dnombrement suivant, sous la forme de la srie $\tilde{D}_3^\bullet(t)$, des sous-groupes du groupe modulaire $\mathrm{PSL}_2(\ZZ)$~:
\begin{align*}
	\tilde{D}_3^\bullet(t) = t+{t}^{2}+4\,{t}^{3}+8\,{t}^{4}+5\,{t}^{5}+22\,{t}^{6}+42\,{t}^{7}
	+\dots
\end{align*}

\begin{table}
\begin{center}
\begin{align*}
\left\{
\begin{aligned}
&\qquad a^*_{0} =1,\quad a^*_{1} =1,\quad a^*_{2} =1,\quad a^*_{3} =2,
\quad a^*_{4} ={\tfrac {15}{4}},\quad a^*_{5} ={\tfrac {91}{20}},
 \\
& \left( {n}^{4}+18\,{n}^{3}+119\,{n}^{2
}+343\,n+366 \right) a^*_{n+6} = \\
& \qquad \left( {n}^{5}+18\,{n}^{4}+121\,{n}
^{3}+373\,{n}^{2}+511\,n+242 \right) a^*_n \\
&\quad + \left( 3\,{n}^{2}+15\,n+18 \right) a^*_{n+1}  \\
&\quad + \left(2\,{n}^{4}+33\,{n}^{3}+205\,{n}^{2}+566\,n+ 582 \right) a^*_{n+2} \\
&\quad + \left( 3\,{n}^{4}+52\,{n}^{3}+938\,n+333\,{n}^{2}+982 \right) a^*_{n+3} \\
&\quad + \left( {n}^{3}+12\,{n}^{2}+53\,n+85 \right) a^*_{n+4} \\
&\quad + \left( {n}^{3}+9\,{n}^{2}+20\,n+1 \right) a^*_{n+5}
\end{aligned}
\right.
\end{align*}
\caption{Relation de rcurrence de longueur six vrifie par les coefficients le la srie $D_3^*(t)$.}
\label{tab:rec:diag:3:non:connexes}
\end{center}
\end{table}

\begin{table}
\begin{center}
\begin{align*}
\tilde{D}_3^\bullet(t)\, &= \,t+{t}^{2}+4\,{t}^{3}+8\,{t}^{4}+5\,{t}^{5}+22\,{t}^{6}+42\,{t}^{7}+40
\,{t}^{8}+120\,{t}^{9}+265\,{t}^{10}+286\,{t}^{11} \\
&\quad+764\,{t}^{12}+1729
\,{t}^{13}+2198\,{t}^{14}+5168\,{t}^{15}+12144\,{t}^{16}+17034\,{t}^{
17}+37702\,{t}^{18} \\
&\quad+88958\,{t}^{19}+136584\,{t}^{20}+288270\,{t}^{21}+
682572\,{t}^{22}+1118996\,{t}^{23} \\
&\quad+2306464\,{t}^{24}+5428800\,{t}^{25}
+9409517\,{t}^{26}+19103988\,{t}^{27}+44701696\,{t}^{28} \\
&\quad+80904113\,{t}
^{29}+163344502\,{t}^{30}+379249288\,{t}^{31}+711598944\,{t}^{32} \\
&\quad+
1434840718\,{t}^{33}+3308997062\,{t}^{34}+6391673638\,{t}^{35}+
12921383032\,{t}^{36} \\
&\quad+29611074174\,{t}^{37}+58602591708\,{t}^{38}+
119001063028\,{t}^{39} \\
&\quad+271331133136\,{t}^{40}+547872065136\,{t}^{41}+
1119204224666\,{t}^{42} \\
&\quad+2541384297716\,{t}^{43}+5219606253184\,{t}^{44
}+10733985041978\,{t}^{45} \\
&\quad+24300914061436\,{t}^{46}+50635071045768\,{t
}^{47}+104875736986272\,{t}^{48} \\
&\quad+236934212877684\,{t}^{49}+
499877970985660\,{t}^{50}+o(t^{50})
\end{align*}
\caption{Dveloppement  l'ordre cinquante de la srie $\tilde{D}_3^\bullet(t)$ donnant le nombre de diagrammes trivalent connexes points (A005133).}
\label{tab:ser:gen;d3dot:labelled}
\end{center}
\end{table}
Il reste maintenant  calculer la srie $\tilde{D}_3(t)$, ce qui sera fait au moyen de sries de Joyal-Plya au paragraphe \ref{calc:expl:prim:meth}. Nous amliorons le calcul au moyen de la notion de forme factorise et nous donnons une formule explicite au paragraphe \ref{calc:expl:deux:meth}.

\begin{application}

Le groupe modulaire admet donc un seul sous-groupe d'indice deux, ncessairement distingu et quatre sous-groupes d'indice trois. Il est facile sur les diagrammes correspondants de constater que parmi ces quatre sous-groupes, un seul est distingu et que les trois autres sont conjugus entre eux.
\begin{center}
\begin{tabular}{ccc}
\begin{tabular}{c}
\includegraphics{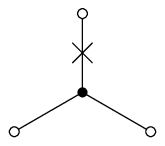}
\end{tabular}
&\hspace{1cm}&
\begin{tabular}{ccccc}
\includegraphics{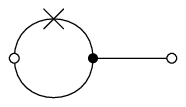}
&\hspace{1cm}&
\includegraphics{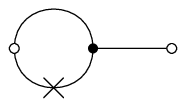}
&\hspace{1cm}&
\includegraphics{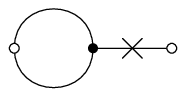}
\end{tabular}
\end{tabular}
\end{center}

On a huit sous-groupes d'indice quatre formant deux classes de conjugaison conrrespondant chacune au quatre faons de pointer chacun des deux diagrammes trivalents de taille quatre. On a cinq sous-groupes d'indice cinq tous conjugus entre eux, correspondant aux cinq pointages de l'unique diagramme trivalent de taille cinq.

Parmi les sous-groupes d'indice six, seuls deux sont distingus, ils correspondent aux deux diagrammes non-quivalents suivants,
\begin{center}
\begin{tabular}{ccc}
\begin{tabular}{c}
\includegraphics{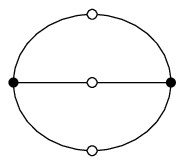}
\end{tabular}
&
\hspace{2cm}
&
\begin{tabular}{c}
\includegraphics{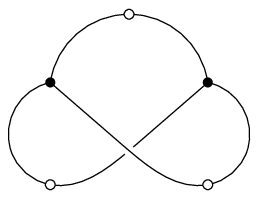}
\end{tabular}
\end{tabular}
\end{center}
qui ont pour groupes d'automorphisme, le groupe symtrique $\mathfrak{S}_3$ et le groupe cyclique $\ZZ/6\ZZ$ respectivement.

Ces deux diagrammes correspondent  deux sous-groupes bien connus de $\mathrm{PSL}_2(\ZZ)$.
Le premier, n'est autre que le \emph{groupe modulaire de niveau deux}, lequel est isomorphe au groupe libre  deux gnrateurs et est associ  la classification des courbes elliptiques modulo les isognies centralises par leurs points de deux-division. On l'obtient comme noyau du morphisme naturel $\mathrm{PSL}_2(\ZZ) \to \mathrm{PSL}_2(\ZZ/2\ZZ) \simeq \mathfrak{S}_3$.
Le deuxime diagramme correspond lui, au sous-groupe driv constitu des commutateurs de $\mathrm{PSL}_2(\ZZ)$, le groupe cyclique $\ZZ/6\ZZ$ est par consquent l'\emph{ablianis} de $\mathrm{PSL}_2(\ZZ)$, \idest son premier groupe d'homologie.
\end{application}

\subsection{Dnombrement de Permutations}
\label{denombrement:permutation}

Nous cherchons  prsent  dterminer explicitement les coefficients des sries $\ZSeries_{S_n}$ pour $n = 2$ et $3$ afin de pouvoir calculer les termes de la srie gnratrice des types de l'espce $D_3$. Nous le faisons plus gnralement pour $n = p$ un nombre premier quelconque, et le mme raisonnement s'appliquerait pour un nombre entier $n$ quelconque avec des formules plus compliques au fur et  mesure que le nombre de diviseurs de $n$ augmente.

\begin{lemme}
On a l'criture suivante pour la srie d'index cyclique $\ZSeries_C$ associe  l'espce des cycles,
\begin{align*}
	\ZSeries_C = \sum_{\substack{r \ge 1 \\ s \ge 1}} \frac{\varphi(r)\,x_r^s}{rs}
\end{align*}
o $\varphi$ dsigne l'indicatrice d'Euler.
\end{lemme}

\begin{demonstration}
A partir de l'isomorphisme naturel $S \simeq \mathrm{Ens}(C)$ exprimant l'existence et l'unicit de la dcomposition des permutations en cycles  support disjoints, on obtient,
au moyen de la formule d'inversion de Mbius,
la srie indicatrice de l'espce des cycles :
\begin{align*}
	\ZSeries_C &= \sum_{k \ge 1} \frac{\mu(k)}{k} \log \left(\prod_{l\ge 1}\frac{1}{1-x_{kl}}\right) \\
		&= \sum_{r \ge 1} \sum_{d \,|\, r} \frac{\mu(d)}{d} \log \left( \frac{1}{1-x_{r}}\right) \\
		&= \sum_{r \ge 1} \frac{\varphi(r)}{r}\log\left(\frac{1}{1-x_r}\right) \\
		&= \sum_{\substack{r \ge 1 \\ s \ge 1}} \frac{\varphi(r)\,x_r^s}{rs}
\end{align*}
Fin de la dmonstration.
\hfill $\Box$ \end{demonstration}

\begin{lemme}
On a l'criture suivante pour la srie d'index cyclique $\ZSeries_{C_n}$ associe  l'espce des cycles de longueur $n$,
\begin{align*}
	\ZSeries_{C_n} = \sum_{rs=n} \frac{\varphi(r)\,x_r^s}{n}
\end{align*}
\end{lemme}

\begin{demonstration}
Il s'agit de regrouper les termes de la srie $\ZSeries_C$ calcule prcdemment, en paquets homognes de mme poids. On obtient alors l'expression de $\ZSeries_{C_n} = \ZSeries_C^{[\,n\,]}$, la composante homogne de poids $n$ de $\ZSeries_C$ :
\begin{align*}
	\ZSeries_C = \sum_{n \ge 0} \sum_{rs=n} \frac{\varphi(r)\,x_r^s}{n}
	\qquad \text{ d'o, } \quad
	\ZSeries_{C_n} = \sum_{rs=n} \frac{\varphi(r)\,x_r^s}{n}
\end{align*}
Fin de la dmonstration.
\hfill $\Box$ \end{demonstration}

\begin{theoreme}
\label{ZSn:factored}
La srie indicatrice de l'espce $S_n$ des permutations $\sigma$ d'ordre $n$ {\rm(}\idest qui satisfont  la relation $\sigma^n=1${\rm)} est donne par la formule ci-dessous o $\varphi$ dsigne l'indicatrice d'Euler.
\begin{align*}
	\ZSeries_{S_n} = \prod_{k\ge 1} \prod_{\substack{rs\,|\,n\\r\,|\,k}} \exp\frac{\varphi(r)\,x_{k}^s}{ks}
\end{align*}
\end{theoreme}

\begin{important}
On observe un phnomne de sparation des variables, lequel est  l'origine de l'algorithme rapide de calcul de sries gnratrices des types d'isomorphisme, donn au paragraphe {\rm \ref{calc:expl:deux:meth}}~{\rm:} la srie gnratrice $\ZSeries_{S_n}$ s'exprime comme produit de sries univaries en les $x_k$.
\end{important}

\begin{demonstration}
Comme une permutation d'ordre $n$ se dcompose canoniquement en cycles d'ordre $d$ o $d$ divise $n$, on obtient un isomorphisme naturel $S_n \simeq \mathrm{Ens} ( \sum_{d\,|\,n} C_d)$. D'o,
\begin{align*}
	\ZSeries_{S_n} &= \exp \sum_{k\ge 1} \frac{1}{k}
		\sum_{d\,|\,n} \sum_{rs = d} \frac{\varphi(r)\,x_{kr}^s}{d} \\
	&= \exp \sum_{k\ge 1} \sum_{rs\,|\,n} \frac{\varphi(r)\,x_{kr}^s}{krs} \\
	&= \prod_{k\ge 1} \prod_{\substack{rs\,|\,n\\r\,|\,k}} \exp\frac{\varphi(r)\,x_{k}^s}{ks}
\end{align*}
Fin de la dmonstration.
\hfill $\Box$ \end{demonstration}

Les deux rsultats suivants donnent les coefficients de la srie $\ZSeries_{S_n}$ sous une forme plus explicite en se restreignant pour simplifier, au cas o l'entier $n$ est un nombre premier $p$. Les raisonnements s'tendent sans difficult au cas d'un entier $n$ quelconque mais nous ne ferons usage que des cas $p=2$ et $p=3$ dans la suite.

\begin{theoreme}
\label{coeff:ZSp:factored}
Pour $p$ un nombre premier quelconque, on a l'criture suivante, sous forme factorise de la srie indicatrice de l'espce des permutations d'ordre $p$,
\begin{align*}
	\ZSeries_{S_p} &= \prod_{k \ge 1}\,\exp \left(\chi_{p,k}\frac{x_k}{k}+\frac{x_k^p}{pk}\right) \\
		&= \prod_{k \ge 1}\left(\sum_{n \ge 0} \;\, \sum_{n_1+pn_2=n} \frac{\chi_{p,k}^{n_1}}{n_1!\,n_2!\,k^{n_1+n_2}\,p^{n_2}} \;\, x_k^n\right)
\end{align*}
o le symbole $\chi_{p,k}$ vaut $p$ ou $1$ suivant que $k$ est divisible par $p$ ou non.
\end{theoreme}

\begin{demonstration}
Il s'agit d'expliciter les coefficients de la srie $\ZSeries_{S_p}$
\begin{equation*}
\begin{split}
	\ZSeries_{S_p} &= \prod_{k\ge 1} \, \exp \sum_{\substack{rs\,|\,p\\r\,|\,k}}\frac{\varphi(r)\,x_{k}^s}{ks} \\
	&= \prod_{k\ge 1} \, \exp \Biggr(
	\frac{\varphi(1)\,x_{k}}{k}
	+ \frac{\varphi(1)\,x_{k}^p}{kp}
	+ \begin{cases}
		\dfrac{\varphi(p)\,x_{k}}{k} & \text{si $k \equiv 0 \mod p$} \\
		0 & \text{sinon}
	\end{cases}
	\Biggl)
	\\
	&= \prod_{k\ge 1} \, \exp \Biggr(
	\chi_{p,k}\frac{x_{k}}{k}
	+ \frac{x_{k}^p}{kp}
	\Biggl)
	\quad\text{ avec $\chi_{p,k} =\begin{cases}
		p & \text{si $k \equiv 0 \mod p$} \\
		1 & \text{sinon.}
	\end{cases}$}
\end{split}
\end{equation*}
Il suffit  prsent d'extraire les coefficients de cette srie indicatrice. On observe pour cela le dveloppement de Taylor suivant,
\begin{align*}
	\exp\,(az+bz^p) = \sum_{\substack{n_1\ge 0 \\ n_2 \ge 0}} \frac{a^{n_1}\,b^{n_2}}{n_1!\,n_2!} z^{n_1+pn_2}
\end{align*}
qui nous donne bien (en effectuant un changement de variables convenable),
\begin{align*}
	\ZSeries_{S_p} 
		&= \prod_{k \ge 1}\left(\sum_{n \ge 0} \;\, \sum_{n_1+pn_2=n} \frac{\chi_{p,k}^{n_1}}{n_1!\,n_2!\,k^{n_1+n_2}\,p^{n_2}} \;\, x_k^n\right)
\end{align*}
Fin de la dmonstration.
\hfill $\Box$ \end{demonstration}

\begin{theoreme}
\label{coeff:ZSp}
Le nombre $u_p(\sigma)$ de permutations d'ordre $p$ premier sur un ensemble  $n$ lments, qui commutent  une permutation $\sigma$ ne dpend que du type cyclique $(\sigma_1, \dots, \sigma_n)$ de $\sigma$ {\rm(}o $\sigma_k$ dsigne le nombre de cycles de longueur $k$ prsent dans la dcomposition cyclique de $\sigma${\rm)}, et vaut~:
\begin{equation*}
	u_{p}(\sigma) = \prod_{k=1}^{n}
	\sum_{n_1+pn_2=\sigma_k} \frac{\sigma_k!\,k^{\sigma_k}\,\chi_{p,k}^{n_1}}{n_1!\,n_2!\,k^{n_1+n_2}\,p^{n_2}}
	\quad\text{avec $\chi_{p,k} = \left\{ \begin{aligned}
		p &\quad \text{si $p$ divise $k$,}  \\
		1 &\quad \text{sinon.}
		\end{aligned}\right.$}
\end{equation*}
\end{theoreme}

\begin{demonstration}
Le premier point est vident. Il est alors permis de noter $u_p(\sigma_1,\dots,\sigma_n)$ la valeur de $u_p(\sigma)$. Pour montrer le deuxime point nous nous appuyons sur le thorme \ref{coeff:ZSp:factored} dmontr prcdemment.
\begin{align*}
	\ZSeries_{S_p} 
		&= \prod_{k \ge 1}\left(\sum_{n \ge 0} \;\, \sum_{n_1+pn_2=n} \frac{\chi_{p,k}^{n_1}}{n_1!\,n_2!\,k^{n_1+n_2}\,p^{n_2}} \;\, x_k^n\right)
\end{align*}
et puisque l'on a, par dfinition~:
\begin{align*}
	\ZSeries_{S_p} = \sum_{n\ge 0}\, \sum_{k_1+2k_2+\dots+nk_n=n}
			\frac{u_p(k_1,\dots,k_n)}{1^{k_1}\,k_1!\cdots n^{k_n}\,k_n!} \, x_1^{k_1}\cdots x_n^{k_n}
\end{align*}
Il ne reste plus qu' identifier les coefficients de ces deux critures.
Fin de la dmonstration.
\hfill $\Box$ \end{demonstration}

\subsection{Calcul explicite : premire mthode}
\label{calc:expl:prim:meth}
Fort des formules que l'on vient d'obtenir nous tchons  prsent de compter les diagrammes trivalents non-points. Nous commenons pour cela par valuer explicitement les sries indicatrices $\ZSeries_{S_2}$ et $\ZSeries_{S_3}$ de l'espce des permutations d'ordre \emph{deux} et \emph{trois} respectivement, puis nous valuons celle de l'espce $D_3^* \simeq S_2 \times S_3$ des diagrammes trivalents non-ncessairement connexes.

Pour donner une ide des calculs, nous prsentons dans les tables numriques \ref{fig:ZS2:coeffs} et \ref{fig:ZS3:coeffs} ci-contre, le dveloppement en filtration \emph{sept} de ces deux sries. Elles sont obtenues en explicitant la formule du thorme \ref{coeff:ZSp} avec $p=2$ puis $p=3$. La table  \ref{fig:ZS2:odot:ZS3:coeffs} donne quant  elle les coefficients de la srie indicatrice des cycles de l'espce $S_2 \times S_3$. Ils s'obtiennent, conformment au lemme \ref{lemme:produit:Hadamard} en effectuant le produit de Hadamard des deux sries prcdentes.

\begin{table}
\begin{center}
\begin{align*}
\begin{split}
\ZSeries_{S_2} &=
1
+ x_{1}
+
\frac{1}{2}
	\left(
2\,x_{{1}}^{2}+2\,x_{{2}}
	\right)
+
\frac{1}{6}
	\left(
4\,x_{{1}}^{3}+6\,x_{{1}}x_{{2}}+2\,x_{{3}}
	\right)
	\\
&\quad+
\frac{1}{24}
	\left(
10\,x_{{1}}^{4}+24\,x_{{1}}^{2}x_{{2}}+8\,x_{{1}}x_{{3}}+18\,x_{{
2}}^{2}+12\,x_{{4}}
	\right)
\\
&\quad+
\frac{1}{120}
	\left(
26\,x_{{1}}^{5}+80\,x_{{1}}^{3}x_{{2}}+40\,x_{{1}}^{2}x_{{3}}+90
\,x_{{1}}x_{{2}}^{2}+60\,x_{{1}}x_{{4}}+40\,x_{{2}}x_{{3}}+24\,x_{{5
}}
	\right)
\\
&\quad+
\frac{1}{720}
	\left(
	\begin{aligned}
&76\,x_{{1}}^{6}+300\,x_{{1}}^{4}x_{{2}}+160\,x_{{1}}^{3}x_{{3}}+
540\,x_{{1}}^{2}x_{{2}}^{2}+360\,x_{{1}}^{2}x_{{4}}\\
&+240\,x_{{1}}
x_{{2}}x_{{3}}+144\,x_{{1}}x_{{5}}+300\,x_{{2}}^{3}+360\,x_{{2}}x_{{
4}}+160\,x_{{3}}^{2}\\
&+240\,x_{{6}}
\end{aligned}
	\right)
\\
&\quad+
\frac{1}{5040}
\left(
\begin{aligned}
&232\,x_{{1}}^{7}+1092\,x_{{1}}^{5}x_{{2}}+700\,x_{{1}}^{4}x_{{3}
}+2520\,x_{{1}}^{3}x_{{2}}^{2}+1680\,x_{{1}}^{3}x_{{4}} \\
&+1680\,x
_{{1}}^{2}x_{{2}}x_{{3}}+1008\,x_{{1}}^{2}x_{{5}}+2100\,x_{{1}}x_{
{2}}^{3}+2520\,x_{{1}}x_{{2}}x_{{4}}\\
&+1120\,x_{{1}}x_{{3}}^{2}+1680
\,x_{{1}}x_{{6}}+1260\,x_{{2}}^{2}x_{{3}}+1008\,x_{{2}}x_{{5}}\\
&+840\,
x_{{3}}x_{{4}}+720\,x_{{7}}
\end{aligned}
	\right)
\\
&\quad+\text{ etc...}
\end{split}
\end{align*}
\caption{Srie indicatrice des cycles de l'espce $S_2$.}
\label{fig:ZS2:coeffs}
\end{center}
\end{table}

\begin{table}
\begin{center}
\begin{align*}
\begin{split}
\ZSeries_{S_3} &=
1
+
x_{{1}}
+
\frac{1}{2}\left(x_{{1}}^{2}+x_{{2}}\right)
+
\frac{1}{6}\left(3\,x_{{1}}^{3}+3\,x_{{1}}x_{{2}}+6\,x_{{3}}\right)
\\
&\quad+
\frac{1}{24}\left(
9\,x_{{1}}^{4}+6\,x_{{1}}^{2}x_{{2}}+24\,x_{{1}}x_{{3}}+3\,x_{{2}}^{2}+6\,x_{{4}}\right)
\\
&\quad+
\frac{1}{120}\left(
\begin{aligned}
21\,x_{{1}}^{5}+30\,x_{{1}}^{3}x_{{2}}+60\,x_{{1}}^{2}x_{{3}}+15\,x_{{1}}x_{{2}}^{2}+30\,x_{{1}}x_{{4}}
+60\,x_{{2}}x_{{3}}+24\,x_{{5}}
\end{aligned}
\right)
\\
&\quad+
\frac{1}{720}\left(
\begin{aligned}
&81\,x_{{1}}^{6}+135\,x_{{1}}^{4}x_{{2}}+360\,x_{{1}}^{3}x_{{3}}+45\,x_{{1}}^{2}x_{{2}}^{2}
+90\,x_{{1}}^{2}x_{{4}}+360\,x_{{1}}x_{{2}}x_{{3}}
\\
&+144\,x_{{1}}x_{{5}}+135\,x_{{2}}^{3}
+90\,x_{{2}}x_{{4}}+360\,x_{{3}}^{2}+360\,x_{{6}}
\end{aligned}
\right)
\\
&\quad+
\frac{1}{5040}\left(
\begin{aligned}
&351\,x_{{1}}^{7}+441\,x_{{1}}^{5}x_{{2}}+1890\,x_{{1}}^{4}x_{{3}}+315\,x_{{1}}^{3}{x_{{2}}}^{2}
+630\,x_{{1}}^{3}x_{{4}}
\\
&+1260\,x_{{1}}^{2}x_{{2}}x_{{3}}+504\,x_{{1}}^{2}x_{{5}}
+945\,x_{{1}}x_{{2}}^{3}+630\,x_{{1}}x_{{2}}x_{{4}}+2520\,x_{{1}}x_{{3}}^{2}
\\
&+2520\,x_{{1}}x_{{6}}
+630\,x_{{2}}^{2}x_{{3}}+504\,x_{{2}}x_{{5}}+1260\,x_{{3}}x_{{4}}+720\,x_{{7}}
\end{aligned}
\right)
\\&\quad+\text{ etc...}
\end{split}
\end{align*}
\caption{Srie indicatrice des cycles de l'espce $S_3$.}
\label{fig:ZS3:coeffs}
\end{center}
\end{table}

\begin{table}
\begin{center}
\begin{align*}
\begin{split}
\ZSeries_{S_2} \odot \ZSeries_{S_3} &=
1
+
x_{{1}}
+
\frac{1}{2}\left(2\,x_{{1}}^{2}+2\,x_{{2}}\right)
+
\frac{1}{6}\left(12\,x_{{1}}^{3}+6\,x_{{1}}x_{{2}}+6\,x_{{3}}\right)
\\
&\quad+
\frac{1}{24}\left(
90\,x_{{1}}^{4}+24\,x_{{1}}^{2}x_{{2}}+24\,x_{{1}}x_{{3}}+18\,x_{{2}}^{2}+12\,x_{{4}}\right)
\\
&\quad+
\frac{1}{120}\left(
\begin{aligned}
&546\,x_{{1}}^{5}+240\,x_{{1}}^{3}x_{{2}}+120\,x_{{1}}^{2}x_{{3}}
+90\,x_{{1}}x_{{2}}^{2}+60\,x_{{1}}x_{{4}} \\
&+120\,x_{{2}}x_{{3}}+24\,x
_{{5}}
\end{aligned}
\right)
\\
&\quad+
\frac{1}{720}\left(
\begin{aligned}
&6156\,x_{{1}}^{6}+2700\,x_{{1}}^{4}x_{{2}}+1440\,x_{{1}}^{3}x_{{
3}}+540\,x_{{1}}^{2}x_{{2}}^{2} \\
&+360\,x_{{1}}^{2}x_{{4}}+720\,x_{
{1}}x_{{2}}x_{{3}}+144\,x_{{1}}x_{{5}}+2700\,x_{{2}}^{3} \\
&+360\,x_{{2}
}x_{{4}}+1440\,x_{{3}}^{2}+720\,x_{{6}}
\end{aligned}
\right)
\\
&\quad+
\frac{1}{5040}\left(
\begin{aligned}
&81432\,x_{{1}}^{7}+22932\,x_{{1}}^{5}x_{{2}}+18900\,x_{{1}}^{4}x
_{{3}}+7560\,x_{{1}}^{3}x_{{2}}^{2} \\
&+5040\,x_{{1}}^{3}x_{{4}}+
5040\,x_{{1}}^{2}x_{{2}}x_{{3}}+1008\,x_{{1}}^{2}x_{{5}}+18900\,x_
{{1}}x_{{2}}^{3} \\
&+2520\,x_{{1}}x_{{2}}x_{{4}}+10080\,x_{{1}}x_{{3}}
^{2}+5040\,x_{{1}}x_{{6}}+3780\,x_{{2}}^{2}x_{{3}} \\
&+1008\,x_{{2}}x_{{
5}}+2520\,x_{{3}}x_{{4}}+720\,x_{{7}}
\end{aligned}
\right)
\\&\quad+\text{ etc...}
\end{split}
\end{align*}
\caption{Srie indicatrice des cycles de l'espce $S_2 \times S_3$.}
\label{fig:ZS2:odot:ZS3:coeffs}
\end{center}
\end{table}

En effectuant  le changement de variable $x_k \to t^k$ dans cette dernire srie on obtient, conformment au lemme de condensation et  l'isomorphisme $D_3^* \simeq S_2 \times S_3$, les huit premiers termes de la srie gnratrice des types d'isomorphisme de l'espce $D_3^*$ des diagrammes trivalents non-ncessairement connexes~:
\begin{align*}
\tilde{D}_3^*(t) = 1+t+2\,{t}^{2}+4\,{t}^{3}+7\,{t}^{4}+10\,{t}^{5}+24\,{t}^{6}+37\,{t}^{7}+\dots
\end{align*}
La srie gnratrice des types d'isomorphismes de l'espce $D_3$ des diagrammes trivalents connexes s'obtient, au moyen de la formule d'inversion de Mbius en vertu de l'isomorphisme naturel $D_3^* \simeq \mathrm{Ens}(D_3)$ qui traduit l'existence et l'unicit de la dcomposition d'un diagramme trivalent quelconque en ses composantes connexes.
\begin{align*}
	\tilde{D}_3(t) = \sum_{k \ge 1} \frac{\mu(k)}{k} \log( \tilde{D}_3^*(t^k))
\end{align*}
Le calcule explicite donne les premiers  coefficients suivants~:
\begin{align*}
\tilde{D}_3(t) = t+{t}^{2}+2\,{t}^{3}+2\,{t}^{4}+{t}^{5}+8\,{t}^{6}+6\,{t}^{7}+\dots
\end{align*}

\subsection{Calcul explicite : deuxime mthode}
\label{calc:expl:deux:meth}

Le nombre et la taille des termes d'une srie indicatrice de cycle sont souvent prohibitifs et conduisent  des calculs difficilement menables. En effets, mme  l'aide de puisants ordinateurs on ne peut en gnral gure esprer calculer en grand poids, car il y a en gnral $p_n$ termes en poids $n$ (o $p_n$ dsigne le nombre de partitions d'un ensemble  $n$ lments) et donc $p_0 + p_1 + p_2 + \dots + p_n$ termes en tout dans une srie tronque en filtration $n$.

A titre indicatif, signalons qu'en filtration cinquante, il y a dj plus d'un million de termes dans une srie d'index cyclique gnrique et que dans les cas habituels, la plupart d'entre eux sont tellement gros que leur criture ncessitent plusieurs pages manuscrites. En poids cinq cents, il faudrait compter environ $4,2 \times 10^{22}$ termes dont la taille serait bien plus grande encore. Autant dire que l'on ne peut gure esprer mener les calculs trs loin sans introduire d'ide nouvelle.

Nous avons russi  mener les calculs au-del du poids cinq cents  l'aide d'un station de bureau en un temps trs bref (environ un quart d'heure). La solution provient de la notion suivante~:

\begin{definition}
\label{def:serie:separable}
Une srie indicatrice de cycle $\ZSeries$ est dite \emph{sparable} si elle admet une criture de la forme suivante :
\begin{equation*}
\label{form:fact:def}
	\ZSeries = \prod_{k \ge 1}\left(\sum_{n \ge 0} \frac{a_{k, n}}{k^n\,n!} \,x_k^n\right)
	\quad \text{avec $a_{k,0} = 1$ pour tout $k\ge 1$.}
\end{equation*}
(un produit de sries univaries en les $x_k$.)
\end{definition}
Une telle criture s'appelle \emph{forme factorise} de la srie $\ZSeries$. Les coefficients $a_{k,n}$ sont alors \emph{uniques}.
En dveloppant, on obtient~:
\begin{equation*}
\begin{aligned}
	\ZSeries &= \sum_{n\ge 0}\, \sum_{k_1+2k_2+\dots+nk_n=n}
			\left(
			\frac{a_{1,k_1}}{1^{k_1}\,k_1!}\, {x_1^{k_1}}
			\right)
			\cdots
			\left(
			\frac{a_{n,k_n}}{n^{k_n}\,k_n!}\, {x_n^{k_n}}
			\right)
			\\
			&= \sum_{n\ge 0}\, \sum_{k_1+2k_2+\dots+nk_n=n}
			\frac{a_{k_1,\dots,k_n}}{1^{k_1}\,k_1!\cdots n^{k_n}\,k_n!}\,x_1^{k_1}\cdots x_n^{k_n}
\end{aligned}
\end{equation*}
o l'on a pos $a_{k_1,\dots,k_n} = a_{1,k_1}\cdots a_{n,k_n}$.

L'intrt principal de cette notion est que le nombre de termes prsents dans chacune des sommes entre parenthses dans la formule de la dfinition \ref{def:serie:separable} est $m/n$ en filtration $m$ puisque le terme $x_n^k$ est de poids $kn$. Le nombre total de termes est donc $m + m/2 + m/3 + \dots + 1 = O(m\,\log(m))$ ce qui est petit. Si l'on est capable d'effectuer les calculs directement sur la forme factorise, on assiste  un effondrement de la complexit, que l'on n'tait pas en mesure d'attendre  priori. C'est  proprement parler, un petit miracle de l'arithmtique.

Une proprit tonnante des sries sparables, non vidente  priori, est prcisment que le produit de Hadamard de deux telles sries s'exprime encore de faon simple en terme des coefficients de leurs formes factorises.

\begin{lemme} De faon prcise, si
\begin{align*}
	\ZSeries_1 = \prod_{k \ge 1}\left(\sum_{n \ge 0} \frac{a_{k, n}}{k^n\,n!} \,x_k^n\right)
	\quad\text{ et }\quad
	\ZSeries_2 = \prod_{k \ge 1}\left(\sum_{n \ge 0} \frac{b_{k, n}}{k^n\,n!} \,x_k^n\right)
\end{align*}
sont deux sries d'index cyclique sparables, prsentes sous leurs formes factorises, alors,
\begin{align*}
	\ZSeries_1 \odot \ZSeries_2 = \prod_{k \ge 1}\left(\sum_{n \ge 0} \frac{a_{k, n}\, b_{k, n}}{k^n\,n!} \,x_k^n\right)
\end{align*}
\end{lemme}
\begin{demonstration}
Immdiat par inspection des formules.
\hfill $\Box$ \end{demonstration}

\begin{application}
On aboutit finalement  la srie gnratrice $\tilde{D}_3(t)$, laquelle numre pour un nombre d'artes donn, les diagrammes trivalents connexes  isomorphisme prs, ou ce qui est quivalent, les classes de conjugaison de sous-groupes d'indice fini correspondants dans le groupe modulaire $\mathrm{PSL}_2(\ZZ)$.
\begin{align*}
	\tilde{D}_3(t) &= \sum_{r \ge 1} \,\frac{\mu(r)}{r}\sum_{k \ge 1}\,\log\left( \sum_{n \ge 0} n!\,k^n\, u_{k,n}\,v_{k, n} \,t^{rkn} \right)
\end{align*}
o les lments $u_{k,n}$ et $v_{k,n}$ sont les coefficients de Taylor de la fonction,
\begin{align*}
	\exp \Biggr(
	\chi_{p,k}\frac{x_{k}}{k}
	+ \frac{x_{k}^p}{kp}
	\Biggl)
	\quad\text{ avec $\chi_{p,k} =\begin{cases}
		p & \text{si $k \equiv 0 \mod p$} \\
		1 & \text{sinon.}
		\end{cases}$}
\end{align*}
pour les valeurs $p=2$ et $p=3$ respectivement. Ce qui donne les valeurs suivantes, en effectuant le dveloppement en srie.
\begin{align*}
	u_{k,n} = \sum_{n_1+2n_2=n} \frac{\chi_{2,k}^{n_1}}{n_1!\,n_2!\,k^{n_1+n_2}\,2^{n_2}} \quad \text{ et } \quad
	v_{k,n} = \sum_{n_1+3n_2=n} \frac{\chi_{3,k}^{n_1}}{n_1!\,n_2!\,k^{n_1+n_2}\,3^{n_2}}
\end{align*}
Lesquelles vrifient des quations de rcurrences videntes, d'ordre deux et trois respectivement.
La table \ref{tab:serie:diag:triv} ci-contre donne les cinquante premiers termes de la srie $\tilde{D}_3(t)$.
\begin{table}
\begin{center}
\begin{equation*}
\begin{split}
\tilde{D}_3(t)\, &= \,t+{t}^{2}+2\,{t}^{3}+2\,{t}^{4}+{t}^{5}+8\,{t}^{6}+6\,{t}^{7}+7\,{t}^{8}+14\,{t}^{9}+27\,{t}^{10}+26\,{t}^{11}\\
&\quad+80\,{t}^{12}
+133\,{t}^{13}+170\,{t}^{14}+348\,{t}^{15}+765\,{t}^{16}+1002\,{t}^{17}+2176\,{t}^{18}\\
&\quad+4682\,{t}^{19}
+6931\,{t}^{20}+13740\,{t}^{21}+31085\,{t}^{22}+48652\,{t}^{23}+96682\,{t}^{24}\\
&\quad+217152\,{t}^{25}
+362779\,{t}^{26}+707590\,{t}^{27}+1597130\,{t}^{28}+2789797\,{t}^{29}\\
&\quad+5449439\,{t}^{30}
+12233848\,{t}^{31}+22245655\,{t}^{32}+43480188\,{t}^{33}\\
&\quad+97330468\,{t}^{34}+182619250\,{t}^{35}+358968639\,{t}^{36}+800299302\,{t}^{37} \\
&\quad+1542254973\,{t}^{38}+3051310056\,{t}^{39}+6783358130\,{t}^{40}+13362733296\,{t}^{41}\\
&\quad+26648120027\,{t}^{42} +59101960412\,{t}^{43}+118628268978\,{t}^{44}\\
&\quad+238533003938\,{t}^{45}+528281671324\,{t}^{46}+1077341937144\,{t}^{47}\\
&\quad+2184915316390\,{t}^{48}+4835392099548\,{t}^{49}+9997568771074\,{t}^{50}
 + o(t^{50})
\end{split}
\end{equation*}
\caption{Dveloppement  l'ordre cinquante de la srie $\tilde{D}_3(t)$ donnant le nombre de diagrammes trivalents connexes (A121350).}
\label{tab:serie:diag:triv}
\end{center}
\end{table}

\end{application}

\begin{table}
\begin{center}
\begin{align*}
[\,t^{500}\,]\,\tilde{D}_3^\bullet(t) = \,
&129430367485890696501112403782149140632007458406669818924 \\
&049655237581302432985235983195547225893573668769081095237 \\
&520334045385563837477539980582454212848418771007253898122 \\
&98261906049050179891685415479424 \\ \\
[\,t^{500}\,]\,\tilde{D}_3(t) = \,
&258860734971781393002224807564298281264014916813339637848 \\
&099310475162604865970471966391094451787371816235545381100 \\
&065419026649727056066352318775170749619149459628751242388 \\
&57108849306258234323621889976 
\end{align*}
\caption{Termes de poids \emph{cinq cents} des sries $\tilde{D}_3^\bullet(t)$ et $\tilde{D}_3(t)$.}
\label{ }
\end{center}
\end{table}

\nocite{*}
\bibliography{art}
\bibliographystyle{plain}

\end{document}